\documentclass[epic,eepic,epsfig,envcountresetchap,envcountsame,draft]{amsart} 
\usepackage{amssymb}
\usepackage{amsmath}
\usepackage{amscd}
\usepackage{epic,eepic,epsfig} 
\makeindex  

\newcommand{\Zb}{{\mathbb Z}}

\def\iii{{\Ic^\prime _2(X,Y)}}
 
\def\Cbb{{ C}}

\def\Nbb{{ N}}
\def\Pbb{{P}}
\def\Qbb{{ Q}}

\def\Zbb{{ Z}}

\def\Cc{{\mathcal C}}

\def\Ic{{\mathcal I}}

\def\Mc{{\mathcal M}}

\def\Uc{{\mathcal U}}
\def\Vc{{\mathcal V}}

\def\tfk{{\frak t}}

\def\01{{\overrightarrow{01}}}
\def\10{{\overrightarrow{10}}}
\def\2{{\overrightarrow{1\infty}}}

\def\hpb{\hfill $\Box$}

\def\lra{\longrightarrow}

\def\al{\alpha}
\def\be{\beta}
\def\de{\delta}
\def\ga{\gamma}
\def\si{\sigma}

\def\Ga{\Gamma}

\def\mod{{\rm mod}}

\def\log{{\rm log}}

\def\pd1{{\partial \Delta [1]}} 
\def\d1{{\Delta [1]}}

\def\zl{\Zbb _\ell}

\def\zlt{\Zbb _\ell^\times}

\def\Cbb{{\mathbb C}}

\def\Nbb{{\mathbb N}}
\def\Pbb{{\mathbb P}}
\def\Qbb{{\mathbb Q}}

\def\Ybb{{\mathbb Y}}
\def\Zbb{{\mathbb Z}}

\def\Cc{{\mathcal C}}

\def\Ic{{\mathcal I}}

\def\Mc{{\mathcal M}}

\def\Uc{{\mathcal U}}
\def\Vc{{\mathcal V}}

\def\Zc{{\mathcal Z}}

\def\z2{{{\Zbb [{\frac{1}{2}}]}}}

\def\dfk{{\mathfrak d}}
\def\ffk{{\mathfrak f}}

\def\tfk{{\mathfrak t}}
\def\zfk{{\mathfrak z}}

\def\al{{\alpha}}
\def\be{{\beta}}
\def\ga{{\gamma}}
\def\si{{\sigma}}
\def\de{{\delta}}
\def\ka{{\kappa}}

\def\De{{\Delta}}

\def\zn{{z^{\frac{1}{\ell ^n}}}}

\def\01{{\overset{\to}{01}}}
\def\10{\overset{\to}{10}}
\def\2{\overset{\to}{1\infty}}

\def\hpb{\hfill $\Box$}

\def\lra{\longrightarrow}

\title{On $\ell$-adic  Galois L-functions}
\author{Zdzis{\l}aw Wojtkowiak}
\begin{document}

\date{\today}

\maketitle

\tableofcontents

\begin{abstract}
Let $z \in \Qbb$ and let $\ga$ be an $\ell$-adic path on  $\Pbb _{\bar \Qbb}^1\setminus \{0,1,\infty \}$ from $\01$ to $z$. For any $\si \in Gal(\bar \Qbb /\Qbb )$, the element
$x^ {-\kappa (\si )}\ffk _\ga (\si )\in \pi _1(\Pbb _{\bar \Qbb}^1\setminus \{0,1,\infty \},\01 )_{pro-\ell}$. 
After the embedding of $\pi _1$ into $\Qbb \{\{X,Y\}\}$ we get the formal power series $\Delta  _\ga (\si )\in  \Qbb \{\{X,Y\}\}$. 
We shall express coefficients of $\Delta  _\ga (\si )$ as integrals over $(\Zbb _\ell )^ r$ with respect to some measures $K_r(z)$. 
The measures $K_r(z)$ are constructed using the tower $\big( \Pbb _{\bar \Qbb}^1\setminus (\{0, \infty \}\cup \mu _{\ell ^ n}\big) _{n\in \Nbb }$ of coverings of
$\Pbb _{\bar \Qbb}^1\setminus \{0,1,\infty \}$. Using the integral formulas we shall show congruence relations between coefficients of the formal power series  $\Delta  _\ga (\si )$.
The congruence relations allow  the construction of $\ell$-adic functions of non-Archimedean analysis, which however rest mysterious. 
Only in the special case of the measures $K_1(\10 )$ and $K_1(-1)$ we recover the familiar Kubota-Leopoldt $\ell$-adic $L$-functions. 
We recover also $\ell$-adic analogues of Hurwitz zeta functions. Hence we get also $\ell$-adic analogues of L-series for Dirichlet characters.

\end{abstract}

\section{Introduction}
\smallskip

\noindent {\bf 0.0 Review of results}  In \cite{W2} we have introduce $\ell$-adic Galois polylogarithms. 
For each $z\in \Qbb$, $l_k(z)$ is a function from $G_\Qbb$ to $\Qbb _\ell$.  
These functions $l_k(z)$ are analogues of the classical polylogarithms $Li_k(z)=\sum _{n=1}^\infty \frac{z^n}{n^k}$. 
In the complex case it is natural to replace $k$ by an arbitrary complex number $s$ and to study a function of two variables $z$ and $s$ defined by the  series $\sum _{n=1}^\infty \frac{z^n}{n^s}$. Notice that for $z=1$ we get the Riemann zeta function $\zeta (s)=\sum _{n=1}^\infty \frac{1}{n^s}$.

\medskip

We would like to replace $k$ in $l_k(z)$ by any $s\in \Zbb _\ell$. We shall be able to do it. However the function we get remains mysterious to us. 
We would like to relate it to an $\ell$-adic non-Archimedean analogue of the complex function $\sum _{n=1}^\infty \frac{z^n}{n^s}$. 
At least we would like to relate its values at positive integers to $\ell$-adic non-Archimedean polylogarithms. We are not able to do this. 
Only in a few special cases we do get the expected results.

\medskip

For $z=\10$ the functions we get, are the Kubota-Leopoldt $\ell$-adic $L$-functions (see \cite{Iw}). The key point is the formula
\begin{equation} \label{eq:l_2k(01)}
 l_{2k}(\10 ) ={\frac{B_{2k}}{2\cdot (2k)!}}(1-\chi ^{2k})  
\end{equation}
proved in \cite{W7}, but stated already in \cite{I1}. In \cite{NW2} there is another proof of the formula \eqref{eq:l_2k(01)}.
We get also familiar functions for $z=-1.$

\medskip

The $\ell$-adic polylogarithm  $l_k(z)$ is by the very definition the coefficient  at $YX^{k-1}$ of the power series
$$
\log \Lambda _\ga \in \Qbb _\ell\{\{X,Y\}\},
$$
where $\ga $ is a path on $\Pbb _{\bar \Qbb}^1\setminus \{0,1,\infty \}$ from $\01 $ to $z$ (see \cite[Definition 11.0.1.]{W2}).  The related function  
\[
 li_k(z)
\]
we define as the coefficient at $YX^{k-1}$ of the power series
\[
 \log \big( \exp (-l(z)_\ga \, X)\cdot \Lambda _\ga\big)\in \Qbb _\ell \{\{X,Y\}\}.
\]

For $z=\10$ and  $\ga$   the canonical path on $\Pbb _{\bar \Qbb }^ 1 \setminus \{ 0,1,\infty \}$ from $\01$ to $\10$, the power series $\Lambda _\ga$ was studied in \cite{D}
and \cite{I}.

In \cite{NW} H. Nakamura and the author have  introduced a certain measure $K_1(z)$ on $\Zbb _\ell$ and   shown that
$$
 li_k(z)=\frac{1}{(k-1)!} \int _{\Zbb _\ell}x^{k-1}dK_1(z).
$$
It has been recovered in this way the Gabber formula of the Heisenberg cover (see \cite{D0}). 

\bigskip

In this paper, for any $r\geq 1$  we construct measures $K_r(z)$ on $(\Zbb _\ell)^r$ which generalize the measure  $K_1(z)$. Then we show that the coefficient at 
\[
X^{a_0}YX^{a_1}YX^{a_2}\ldots X^{a_{r-1}}YX^{a_r}  
\]
of the power series 
$$
\log \big( \exp (-l(z)_\ga \, X)\cdot \Lambda _\ga\big)\in \Qbb _\ell \{\{X,Y\}\} 
$$ 
is given by the integral
\begin{equation}\label{eq:start}
{\frac{1}{a_0!a_1!\ldots a_r!}}\int _{(\Zbb _\ell)^r}(-x_1)^{a_0}(x_1-x_2)^{a_1} \ldots (x_{r-1}-x_r)^{a_{r-1}}(x_r)^{a_r}dK_r(z)\, . 
\end{equation}

Using this integral expression we shall be able to prove congruence relations between coefficients of the power series 
$\log \big( \exp (-l(z)_\ga \, X)\cdot \Lambda _\ga\big)$.

In the integral \eqref{eq:start}, after some modifications, we can replace the integers $a_0,\ldots ,a_r$ by arbitrary $s_0,\ldots ,s_r$ in $\Zbb _\ell$. 
However the obtained functions are mysterious. As we already mentioned, only for $r=1$ and $z=\10$ we do get the familiar Kubota-Leopoldt $\ell$-adic L-functions.
The familiar functions we get also for $r=1$ and $z=-1$.
 
Let $\xi _m=e^{\frac{2 \pi i}{m}}$. We assume that $\ell$ does not divide $m$. 
Then using measures $K_1(\xi _m^{-k})\pm K_1(\xi _m^{k})$ we get $\ell$-adic analogues of Hurwitz zeta function. 
Hence we get also $\ell$-adic analogues of L-series for Dirichlet characters.
 
Below we fix notations and conventions used in the paper. We review also the definitions of $\ell$-adic polylogarithms and measures.

\bigskip

\noindent
{\bf 0.1  Notations  and conventions} 
Throughout the paper we fix the following notation and conventions.

We fix a rational prime $\ell$. If $V$ is an algebraic variety over a number field $K$ and $v$ and $z$ are $K$-points or tangential points defined over $K$ we denote by 
\[
 \pi _1(V_{\bar K},v)
\]
the maximal pro-$\ell$ quotient of the \'etale fundamental group of $V_{\bar K}$ based at $v$ and by
\[
 \pi (V_{\bar K};z,v)
\]
the   $\pi _1(V_{\bar K},v)$-torsor of $\ell$-adic paths on $V_{\bar K}$ from $v$ to $z$. We recall that an $\ell$-adic path $\ga$ from $v$ to $z$ on $V_{\bar K}$ 
is an isomorphism of fiber functors $\ga :F_v\to F_z$.

If $\al$ is an $\ell$-adic path from $a$ to $b$ and $\beta $ from $b$ to $c$ then
\[
 \beta \cdot \alpha
\]
is an $\ell$-adic path from $a$ to $c$.

When we  speak about a multiplicative embedding $E$ of $\pi _1$ into an algebra of formal power series we mean that 
\[
 E(\be \cdot \al )=E(\be )\cdot E(\al ) .
\]

We assume that $\bar K\subset \Cbb$. Then we have the comparison homomorphism
\[
 \pi _1(V(\Cbb ),v)\to \pi _1(V_{\bar K},v)
\]
and the comparison map
\[
 \pi  (V(\Cbb );z,v)\to \pi  (V_{\bar K};z,v)\, .
\]

\smallskip
 
In this paper path, homotopy class of path and $\ell$-adic path mean exactly the same. They mean an $\ell$-adic path as defined above. We usually shall say path if we 
can take an element of $\pi _1(V(\Cbb ),v)$ or $\pi  (V(\Cbb );z,v)$.

\smallskip

If $\si \in G_K$ and $\ga$ is a path then
\[
 \si (\ga )=\si \circ \ga \circ \si ^{-1}\, .  
\]
The action of $\pi _1$ and $G_K$ on germs of algebraic functions is the left action.

We define
\[
 \ffk _\ga (\si ):=\ga ^{-1}\cdot \si (\ga )\in \pi _1(V_{\bar K},v)\, .
\]

We denote by
\[
 \Nbb
\]
the set of positive integers and $0$.
For $\al \in \Qbb _\ell$ and $k\in \Nbb$ we denote by 
\[
 C_k^\al
\]
the binomial coefficients. For any positive integer $m$ we set
\[
 \xi _{m}:=e^{\frac{2\pi \sqrt{-1}}{m}}\, .
\]

\medskip

\noindent
{\bf  0.2 Algebraic preliminaries and $\ell$-adic polylogarithms }  We denote by 
\[
 \Qbb _\ell \{\{X,Y\}\}
 \]                                                                                
the $\Qbb _\ell$-algebra of formal power series in two non-commuting variables $X$ and $Y$. The set of Lie  polynomials in
 $\Qbb _\ell \{\{X,Y\}\}$ we denote by $Lie (X,Y)$. It is a free Lie algebra on $X$ and $Y$. The set of formal Lie power series in $\Qbb _\ell\{\{X,Y\}\}$ 
we denote by $L(X,Y)$. The vector space $L(X,Y)$ is a Lie algebra, the completion of $Lie (X,Y)$ with respect to the filtration given by the lower central series.
We denote by
\[I_2\] the closed Lie ideal of $L(X,Y)$ generated by Lie brackets with two or more $Y$'s.

Let $A,B$ be elements of a Lie algebra. We shall use the following inductively defined short hand notation
\[
 [B,A^{(0)}]:=B\;\;{\rm and}\;\;[B,A^{(n+1)}]:=[[B,A^{(n)}],A]\;\;{\rm if }\;\;n\geq 0.
\]
If $P$ is a formal power series without a constant term we shall write $\exp P$ or $e^P$ to denote the formal power series
\[
 \sum _{n=0}^\infty {\frac{P^n}{n!}}\,.
\]

Let $A,B\in L(X,Y)$. The formula
\[
 A\bigcirc B:=\log (\exp A \cdot \exp B)
\]
defines a group multiplication in the set $L(X,Y)$ and it is called the Baker-Campbell
-Hausdorff product. In the group $L(X,Y)$ one has
\[
 A\bigcirc (-A)=0\;.
\]
If $\al \in \Qbb _\ell$ then one can raise elements of the group $L(X,Y)$ to the power $\al$ and
\[
 A^\al =\al A\;.
\]
We denote by 
\[
 \iii
\]
the closed ideal of $ \Qbb _\ell \{\{X,Y\}\}$ generated by all monomials with two $Y$'s and by monomials $X^iY$ for $i>0$.

\medskip

The well known formulas
\[
 X\bigcirc Y\equiv X+Y  {\frac{X}{\exp X-1}}\;\;{\rm mod}\;\;\iii 
\]
and
\[
 Y\bigcirc X\equiv X+Y   {\frac{X\exp X}{\exp X-1}}\;\;{\rm mod}\;\;\iii
\]
are easy consequences of the next lemma.

\medskip

\noindent{\bf Lemma 0.2.1.}  Let $\al, \be \in \Qbb _\ell ^\times$ and let $A$ and $B$ belong to $L(X,Y)$. We assume that
\[
 A\equiv \al X+Y  \Phi _1(X)\;\;{\rm mod}\;\;\iii \;\;{\rm and}\;\;B\equiv \be X+Y  \Phi _2(X)\;\;{\rm mod}\;\;\iii\,,
\]
where $\Phi _1(X)$ and $\Phi_2(X)$ are power series in $X$. Then we have
\[
 A\bigcirc B\equiv Y  \big(\Phi _1(X)  {\frac{\exp (\al X)-1}{\al X}}  e^{\be X}+\Phi _2(X) 
{\frac{\exp (\be X)-1}{\be X}}\big)\cdot
\]
$$
  {\frac{(\al +\be)X}{\exp ((\al +\be)X)-1}}\,+(\al +\be )X\;\;{\rm mod}\;\;\iii\,.
$$
( If the constant $\ga =0$  then the  power series $ {\frac{\exp (\ga X)-1}{\ga X}}$ is equal $1$.)

\medskip

\noindent{\bf Proof.}  We omit the proof of the lemma, which is 
the standard calculation on formal power series. It is similar 
to the proof of the two well known formulas given above. \hpb

\medskip

In the Lie algebra $ L(X,Y)$ we set 
\[
 Z:=-\log (e^X  e^Y)\,.
\]
Then $Z\equiv -X-Y {\frac{X}{\exp X-1}}$ modulo $\iii$.

\bigskip

We recall the definition of $\ell$-adic polylogarithms (see \cite{W2}). 
Let $x$ and $y$ be the generators of the free pro-$\ell$ group $\pi _1(\Pbb ^1_{\bar \Qbb}\setminus \{ 0,1,\infty \},\01 )$ as on   Picture 1.

\[
 \;
\]

$$
\;
$$

\[
\;
\]

$${\rm Picture \;1}$$
Let
\[
 E:\pi _1(\Pbb ^1_{\bar \Qbb}\setminus \{ 0,1,\infty \},\01 )\to \Qbb _\ell \{\{X,Y\}\}
\]
be the continuous multiplicative embedding defined by
\[
 E(x)=\exp X\;\;\;{\rm and}\;\;\;  E(y)=\exp Y\, .
\]
Let $z$ be a $\Qbb$-point or a tangential point defined over $\Qbb$ of $\Pbb ^1 \setminus \{ 0,1,\infty \}$. 
Let $\ga $ be an $\ell$-adic path from $\01$ to $z$ on $\Pbb ^1_{\bar \Qbb}\setminus \{ 0,1,\infty \}$ and let $\si \in G_{\Qbb}$. We set
\[
 \Lambda _\ga (\si ) :=E(\ffk _\ga (\si ))\in \Qbb _\ell \{\{X,Y\}\}.
\]
The formal power series $\log \Lambda _\ga (\si )$ is a Lie series. We defined $\ell$-adic Galois polylogarithms $l_n(z)_\ga :G_{\Qbb}\to \Qbb _\ell$ by the congruence
\begin{equation}\label{eq:defpoly}
 \log \Lambda _\ga (\si )\equiv l(z)_\ga (\si )X+\sum _{n=1}^\infty l_n(z)_\ga (\si )[Y,X^{(n-1)}]\;\; {\rm mod} \;\; I_2\, .
\end{equation}
The $\ell$-adic logarithm $l(z)_\ga$ is the Kummer character $\kappa (z)$ associated to $z$ and $l_1(z)_\ga=\kappa (1-z)$.
 
\medskip

Another version of $\ell$-adic polylogarithms 
\[
 li_n(z)_\ga :G_\Qbb \to \Qbb _\ell
\]
we define by the congruence
\[
 \log \big( \exp (-l(z)_\ga (\si )X)\cdot \Lambda _\ga (\si )\big) \equiv \sum _{n=1}^\infty li_n(z)_\ga (\si )[Y,X^{(n-1)}]\;\;{\rm mod }\;\;I_2.
\]

\medskip

The relation between these two versions of $\ell$-adic polylogarithms is given by the equality of formal power series
\begin{equation}\label{eq:litol}
 \sum _{n=1}^\infty li_n(z)_\ga X^{n-1}=\big( \sum _{n=1}^\infty l_k(z)_\ga X^{n-1}\big)  {\frac{\exp (l(z)_\ga X)-1}{l(z)_\ga X}}\;,
\end{equation}
which  follows from Lemma 0.2.1.

The functions 
\[
t_i(z)_\ga :G_\Qbb \to \Zbb _\ell 
 \]
are defined by the congruence
\begin{equation} \label{eq:congruencefort}
 x^{-l(z)_\ga (\si )}\cdot \ffk _\ga (\si )\equiv \prod _{i=1}^\infty (y,x^{(i-1)})^{t_i(z)_\ga (\si )}
\end{equation}
modulo commutators with two or more $y$'s 
 and where \[(y,x):=yxy^{-1}x^{-1}\, ,\;\;  (y,x^{(0)}):=y\;\;{\rm and} \;\;(y,x^{(i+1)}):=((y,x^{(i)}),x)\]
 for $i\geq 1$ (see also \cite{W6}, where these exponents are studied).

\medskip

\noindent
{\bf  0.3 Measures}   In this subsection we collect some elementary properties of measures. Let $X$ be a projective limite of finite sets equipped with 
the limit topology. Further we shall call such  $X$ a profinite set. We denote by 
\[
 CO(X)
\]
the set of compact-open subsets of $X$. A measure $\mu$ on $X$ is a bounded finitely additive function
\[
 \mu :CO(X)\to \Qbb _\ell.
\]

\bigskip

Let $X$ and $Y$ be profinite sets and let $\phi :X\to Y$ be a continuous map. Let $\mu$ be a measure on $X$. We define a measure
\[
 \phi _!(\mu ):CO(Y)\to \Qbb _\ell
\]
on $Y$ by
\[
 (\phi_!\mu )(\Uc ):=\mu (\phi ^{ -1}(\Uc )).
\]
For any $f\in \Cc (Y,\Qbb _\ell)$ -- $\Qbb _\ell$-vector space of continuous functions from $Y$ to $\Qbb _\ell$ -- we have
\begin{equation}\label{eq:fcirc phi}
 \int _Yfd(\phi_!\mu )=\int _X(f\circ \phi)d\mu .
\end{equation}

\bigskip

Let $X$ and $Y$ be profinite sets and let $\phi :X\to Y$ be a continuous open injective map. Let $\nu$ be a measure on $Y$. We define a measure 
\[
 \phi ^!\nu :CO(X)\to \Qbb_\ell
\]
on $X$ by
\[
 (\phi ^!\nu )(\Vc ):=\nu (\phi (\Vc )).
\]
For any $f\in \Cc (Y,\Qbb _\ell )$ we have
\begin{equation}\label{eq:fcirc phi2}
 \int _X(f\circ \phi )d (\phi ^!\nu )=\int _Y (\chi _{\phi (X)}  f)d\nu ,
\end{equation}
where $\chi _A$ is the characteristic function of a subset $A$.

\smallskip

If $\phi$ is a homeomorphism then
\[
 \phi ^!\nu =(\phi ^{-1})_!\nu .
\]
Let $\Uc$ be a compact-open subset of $Y$. Let $i:\Uc \to Y$ be the inclusion. Then the measure $i^!\nu$ we denote also by $\nu _{ \mid \Uc }$.
For $f\in \Cc (Y,\Qbb _\ell)$ we have
\[
 \int _\Uc (f\circ i)d( \nu _{ \mid \Uc } )=\int _Y(\chi _\Uc  f)d\nu .
\]

\medskip

For the profinite set
\[
 X=(\Zbb _\ell )^r
\]
we shall review several equivalent definitions of measure.

\medskip

\noindent
{\bf Definition 0.3.1. } A measure $\mu$ on $(\Zbb _\ell )^r$ is a family of functions
\[
 \big( \mu ^{(n)}:(\Zbb  /\ell ^n \Zbb )^r\to \Qbb _\ell \big) _{n\in \Nbb }
\]
satisfying the distribution relations and which are uniformly bounded. 

\medskip

Therefore the values of all functions $\mu ^{(n)}$ are in ${\frac{1}{\ell ^N}}\Zbb _\ell$ for some $N\geq 0$. For simplicity we shall assume farther that these values are
in $\Zbb _\ell$.

Observe that
\[
 \big( \sum _{\iota \in (\Zbb /\ell ^n)^r}\mu ^{(n)}(\iota ) \iota \big) _{n\in\Nbb }\in {\varprojlim} _n \Zbb _\ell [(\Zbb /\ell^n\Zbb )^r]=\Zbb _\ell [[(\Zbb _\ell )^r]].
\]
Hence we have the following definition.

\medskip

\noindent
{\bf Definition 0.3.2. } A measure $\mu$ on $(\Zbb _\ell )^r$ is an element 
\[
 \mu \in \Zbb _\ell [[(\Zbb _\ell )^r]] .
\]

\medskip

The Iwasawa algebra $\Zbb _\ell[[(\Zbb _\ell)^r]]$ is isomorphic to the algebra of commutative formal power series $\Zbb _\ell[[A_1,A_2\ldots A_r]]$.  
The isomorphism of $\Zbb _\ell$-algebras
\[
 P:\Zbb _\ell[[(\Zbb _\ell)^r]] \to \Zbb _\ell[[A_1,A_2\ldots A_r]]
\]
is given by
\[
 P\big( (\al _1,\al _2\ldots \al _r) \big)=\prod _{i=1}^r(1+A_i)^{\al _i},
\]
for $ (\al _1,\al _2\ldots \al _r) \in (\Zbb _\ell )^r$ and is extended by continuity.
If $\mu \in \Zbb _\ell[[(\Zbb _\ell)^r]]$ then 
\begin{equation} \label{eq:P(mu)}
 P(\mu )(A_1,\ldots ,A_r)= 
\end{equation}
\[
\sum _{n_1=0}^\infty \ldots \sum _{n_r=0}^\infty  \big( \int _{ (\Zbb _\ell)^r }C_{n_1}^{x_1}  C_{n_2}^{x_2} \ldots C_{n_r}^{x_r}
d\mu (x_1,\ldots ,x_r)\big) A_1^{n_1}   A_2^{n_2}\ldots A_r^{n_r} \;.
\]

Let 
\[
 F:\Zbb _\ell[[(\Zbb _\ell)^r]] \to \Qbb _\ell[[X_1,X_2\ldots X_r]]
\]
be given by 
\[
 F(\mu )(X_1,\ldots ,X_r):=P(\mu )(\exp (X_1)-1,\ldots ,\exp (X_r)-1)\; .
\]

Then we have
\begin{equation}\label{eq:F(mu)}
 F(\mu )(X_1,\ldots ,X_r)=
\end{equation}
\[
\sum _{n_1=0}^\infty\ldots \sum _{n_r=0}^\infty {\frac{1}{n_1!  n_2!\ldots n_r!}}\big( \int _{  (\Zbb _\ell)^r }  x_1^{n_1}  
 x_2^{n_2}\ldots  x_r^{n_r}d\mu (x_1,\ldots ,x_r)\big) X_1^{n_1}  X_2^{n_2}\ldots X_r^{n_r}\; .
\]
(see also \cite[pages 290 and 291]{NW})
\medskip

Let 
\[
 \phi :(\Zbb _\ell )^r\to (\Zbb _\ell )^r
\]
be a morphism of $\Zbb _\ell$-modules. We denote by
\[
 \phi ^{(n)} :(\Zbb /\ell ^n\Zbb )^r\to (\Zbb /\ell ^n\Zbb )^r
\]
the induced morphism.  The morphisms  $\phi ^{(n)}$ induce morphisms of group rings
\[
( \phi ^{(n)}) _*:\Zbb _\ell[(\Zbb /\ell ^n\Zbb )^r ] \to \Zbb _\ell[(\Zbb /\ell ^n\Zbb )^r ] 
\]
and in consequence the morphism of the Iwasawa algebras
\[
 \phi _*:\Zbb _\ell[[(\Zbb _\ell)^r]] \to \Zbb _\ell[[(\Zbb _\ell)^r]]\,.
\]

\medskip

\noindent
{\bf Proposition 0.3.3.}   Let $\mu $ be a measure on $(\Zbb _\ell  )^r$. Then we have
\[
 \phi _!(\mu )=\phi _*(\mu )\, .
\]

\medskip

\noindent {\bf Proof.} The element 
\[
 \phi _*(\mu )= \big( (\phi  ^{(n)}) _*)(\mu )\big) _{n\in \Nbb}\in {\varprojlim} _n \Zbb _\ell[(\Zbb /\ell ^n\Zbb )^r ]\, .
\]
We have
\[
 (\phi  ^{(n)})_*(\mu )= (\phi ^{(n)})  _*(\mu ^ {(n)})= (\phi ^{(n)} ) _*(\sum _{\iota \in (\Zbb /\ell ^n\Zbb )^r}\mu ^ {(n)}(\iota )  \iota )=
\sum  _{\iota \in (\Zbb /\ell ^n\Zbb )^r}\mu ^{(n)}(\iota )\phi ^{(n)}( \iota)
\]
\[
 =\sum  _{\kappa  \in (\Zbb /\ell ^n\Zbb )^r} \big(\sum _{\iota \in (\phi ^{(n)} )^{-1}(\kappa )}\mu ^{(n)} (\iota )\big)\kappa\, .
\]

Let $0\leq k_1,\ldots ,k_r<\ell ^n$. Therefore we get
\[
 (\phi _*\mu )\big( (k_1,\ldots ,k_r)+\ell  ^n(\Zbb _\ell )^r\big)=\big( (\phi  ^{(n)}) _*)(\mu )\big) (k_1,\ldots ,k_r)=
\] 
\[ 
\sum _{\iota \in  (\phi ^{(n)})^{-1}(k_1,\ldots ,k_r)}\mu ^{(n)}(\iota )=\mu (\phi ^{-1}\big( (k_1,\ldots ,k_r)+\ell  ^n(\Zbb _\ell )^r\big) )=
\]
\[
(\phi _!\mu )\big( (k_1,\ldots ,k_r)+\ell  ^n(\Zbb _\ell )^r\big) \,.
\]
\hpb

\medskip

\noindent
{\bf Corollary  0.3.4.} Let $A=(a_{i,j})$ be the matrix of $\phi :(\Zbb _\ell )^r\to (\Zbb _\ell )^r$. Then we have
\[
 P(\phi _!\mu )(A_1,\ldots ,A_r)=P(\mu )\big( \prod_{i=1}^r(1+A_i)^{a_{i,1}},\ldots , \prod_{i=1}^r(1+A_i)^{a_{i,r}}\big)
\]
and
\[
  F(\phi _!\mu )(X_1,\ldots ,X_r)=F(\mu )\big( \sum _{i=1}^ra_{i1}X_i,\ldots ,\sum _{i=1}^ra_{ir}X_i\big)\, .
\]

\medskip

If $q\in \Zbb _\ell $ then $\langle q\rangle$ is a positive integer such that $0\leq \langle q\rangle <\ell ^n$ and $\langle q\rangle\equiv q $ modulo $\ell ^n$.

Below we give an example of a measure on $\Zbb _\ell$ which will frequently appear in this paper.

\noindent
{\bf Example 0.3.5.} Let $c\in \Zbb _\ell ^\times $. The Bernoulli measure 
\[
 E_{1,c}=\big(E_{1,c}^{(n)} :\Zbb /\ell ^n\Zbb \to \Qbb _\ell \big)_{n\in \Nbb }
\]
on $\Zbb _\ell$ is defined by
\[
 E_{1,c}^{(n)} (i)={\frac{i}{\ell ^n}}-c{\frac{\langle c^{-1}  i\rangle }{\ell ^n}}+{\frac{c-1}{2}}
\]
for $0\leq i <\ell ^n$.

\medskip
 

\section{Action of the absolut Galois group on  fundamental groups}

\smallskip
 Let $V:=\Pbb ^1_{\bar \Qbb }\setminus (\{0,\infty \}\cup \mu _{\ell ^n})$. We recall that $\ell$ is a fixed prime and that $\pi _1(V ,\01 )$ is the maximal pro-$\ell$ quotient of the \'etale fundamental group of $V $ based at $\01 $. 
We describe the Galois action on generators of $\pi _1(V ,\01 )$. In contrast with our other papers (\cite{W3}, \cite{W4}), we are studying the action of $G_\Qbb$, not merely of $G_{\Qbb (\mu _{\ell ^n})}$. First we recall the construction of generators of $\pi _1(V ,\01 )$.

$$
\;
$$

\[
 \;
\]

\[
 \;
\]

$$\rm Picture \;2$$

\medskip

Let $x\in \pi _1(V ,\01 )$, $y^\prime _k\in \pi _1(V,{\overset{\to}{\xi _{\ell ^n}^k0}} )$ and let  $\beta _k$ be a path from 
$\01$ to ${\overset{\to}{\xi _{\ell ^n}^k0}}$  as on the picture. Let us set 
$$
y_k:=\beta _k^{-1}\cdot y_k^\prime \cdot \beta _k.
$$
Then 
\[
 x,\; y_0,\;y_1,\ldots ,\;y_{\ell^n-1}
\]
are free generators of $\pi _1(V ,\01 )$.
\bigskip

\noindent{\bf Theorem 1.1.}  The Galois group $G_{  \Qbb}$ acts on  $\pi _1(V,\01 )$. For any $\si \in G_\Qbb$ we have
\[
 \si (x)=x^{\chi (\si )}
\]
and
$$
\si (y_k)=((\beta  _{k\cdot \chi (\si )})^{-1}\cdot \si (\beta _k))^{-1}\cdot (y _{k\cdot \chi (\si)})^{\chi (\si )}\cdot ((\beta  _{k\cdot \chi (\si )})^{-1}\cdot \si (\beta _k))
$$
for $k=0,1,\ldots ,\ell^n-1$.

\medskip

\noindent{\bf Proof.}  The Galois group $G_\Qbb$ permutes the missing points $\{0,\infty \}\cup \mu _{\ell ^n}$. Hence it follows that $G_\Qbb$ acts on 
$\pi _1(V ,\01 )$. Let $z$ be the standard coordinate on $\Pbb ^1$. Then $\si \cdot y_k^\prime \cdot \si ^{-1}$ transforms 
$(1-\xi ^{-k\chi (\si )}_{\ell ^n}z)^{\frac{1}{\ell^m}}$ to 
$(1-\xi ^{-k }_{\ell ^n}z)^{\frac{1}{\ell^m}}$, next to $\xi ^1_{\ell^m}(1-\xi ^{-k }_{\ell ^n}z)^{\frac{1}{\ell^m}}$ and finally to 
$\xi _{\ell^m}^{\chi (\si)}(1-\xi ^{-k\chi (\si )}_{\ell ^n}z)^{\frac{1}{\ell^m}}$. Hence it follows that $\si (y_k^\prime )=(y_{k\chi (\si)}^\prime )^{\chi (\si)}$.

 We have 
\[
\si (y_k)=\si (\be _k^{-1}\cdot y_k^\prime \cdot \beta _k)=\si (\be _k^{-1})\cdot \si (y_k^\prime )\cdot \si (\beta _k )=
\]
\[
(\si (\beta _k)^{-1}\cdot \beta _{k\cdot \chi (\si )})\cdot (\beta _{k\cdot \chi (\si )})^{-1}\cdot \si (y_k^\prime )\cdot \beta _{k\cdot \chi (\si )}\cdot 
((\beta _{k\cdot \chi (\si )})^{-1}\cdot \si (\beta _k))=
\] 
\[ 
\big( (\beta  _{k\cdot \chi (\si )})^{-1}\cdot \si (\beta _k)\big) ^{-1}\cdot ( y _{k\cdot \chi (\si)})^{\chi (\si )}\cdot \big( (\beta  _{k\cdot \chi (\si )})^{-1}\cdot \si (\beta _k)\big).
\]   \hpb

\bigskip


\section{Measures associated to towers of projective lines}
\smallskip
In this section we construct measures on $(\Zbb _\ell )^r$, which generalize the measure constructed in \cite{NW}. Next we generalize the principal result of \cite{NW} expressing the $\ell$-adic polylogarithms $li_k(z)$ as the integrals over $\Zbb _\ell$.

\bigskip

\noindent
For each $n\geq 0$ we set
$$
V_n:=\Pbb ^1_{\bar \Qbb }\setminus (\{0,\infty \}\cup \mu _{\ell ^n}).
$$
Let \[f_n^{m+n}:V_{m+n}\to V_n\] be given by \[f_n^{m+n}(z)=z^{\ell^m}.\] Observe that $f_n^{m+n}(\01)=\01$. Hence we get a family of homomorphisms
\begin{equation} \label{eq:comp homo}
(f_n^{m+n})_*:\pi _1(V_{m+n},\01)\to \pi _1(V_{ n},\01)
\end{equation}
satisfying $$(f^{m+n+p}_p)_*=(f^{n+p}_p)_*\circ (f^{m+n+p}_{n+p})_*.$$

Observe that the Galois group $G_\Qbb $ acts on each $\pi _1(V_n,\01)$ and that $(f_n^{m+n})_*$ are $G_\Qbb $-maps. We choose generators 
$$
x_n,\; y_{n,0},\; y_{n,1},\ldots ,\; y_{n,\ell ^n-1}
$$ 
of  
$\pi _1(V_n,\01)$ as in Section 1, i.e. $x_n=x$ and $y_{n,i}=y_i$ in the notation of Section 1. Then we have
\begin{equation} \label{eq:f on gener}
(f_n^{m+n})_*(x_{m+n})=(x_n)^{\ell ^m}\;\;{\rm and}\; \;(f_n^{m+n})_*(y_{m+n,k})=x^{-g}\cdot y_{n,k^\prime}\cdot x^g,  
\end{equation}
where $k=k^\prime +g  \ell^n$ and $0\leq k^\prime <\ell^n$.

Let us set
\[\Ybb _n:=\{X_n,Y_{n,i}\;\mid \; 0\leq i<\ell ^n\}\] and let \[ \Qbb _\ell\{\{\Ybb _n\}\}\] be a $\Qbb _\ell$-algebra of formal power series in non-commuting 
variables $$X_n,\;Y_{n,0},\;Y_{n,1},\ldots ,\; Y_{n,\ell^n-1}.$$ Let
$$
E_n:\pi _1(V_n,\01)\to \Qbb _\ell\{\{\Ybb _n\}\}
$$
be a continuous multiplicative embedding given by \[ E_n(x_n):=\exp X_n\;\;{\rm and}\;\;  E_n(y_{n,i}):=\exp Y_{n,i}\;\;{\rm for}\;\; 0\leq i<\ell ^n.\]
The action of $G_\Qbb $ on $\pi _1(V_n,\01)$ induces the action of $G_\Qbb$ on $\Qbb _\ell\{\{\Ybb _n\}\}$. The
homomorphisms \eqref{eq:comp homo} induce
$G_\Qbb$-morphisms
$$
(f_n^{m+n})_*:\Qbb _\ell\{\{\Ybb _{m+n}\}\}\to \Qbb _\ell\{\{\Ybb _n\}\}
$$ 
such that
$$
(f_n^{m+n})_*\circ E_{m+n}=E_n\circ (f_n^{m+n})_*
$$
and
$$
(f_p^{m+n+p})_*=(f_p^{n+p})_*\circ (f_{n+p}^{m+n+p})_*.
$$
It follows from \eqref{eq:f on gener} that  
\begin{equation} \label{eq:fonYn}
(f_n^{m+n})_*(X_{m+n})=\ell^mX_n\;\;{\rm and}\;\;(f_n^{m+n})_*(Y_{m+n,k})=\exp (-gX)\cdot Y_{n,k^\prime}\cdot \exp (gX),
\end{equation}
if $k=k^\prime +g  \ell^n$ and $0\leq k^\prime <\ell^n$.

\bigskip

\noindent
Let $\al \in \Zbb _\ell$. Then $\al =\sum _{i=0}^\infty \al _i\ell^i$ where $0\leq \al _i<\ell $.
We define \[\al (n):=\sum _{i=0} ^{n-1}\al _i \ell^i.\]

Observe that  $\xi _{\ell ^n}^\al$ is well defined and $(\xi _{\ell^{m+n}}^\al )^{\ell^m}=\xi _{\ell ^n}^\al $. 
Let $g_\al ^{(n)}:V_n\to V_n$ be given by $g_\al ^{(n)}(z)=\xi _{\ell ^n}^\al   z$.

Let $0\leq q <\ell^n$. Let $s_q$ be a path on $V_n$ from $\01$ to ${\overset{\to}{0\xi _{\ell ^n}^q}}$ as on the picture.

$$
\;
$$

\[
 \;
\]
 
\[
 \;
\]

$$
{\rm Picture\;3}
$$

  We define
$$
(x_n)^{{\frac{1}{\ell ^n}}\al}:=s_{\al (n)}\cdot (x_n)^{{\frac{1}{\ell ^n}}(\al -\al (n))}.
$$

\medskip

\noindent
Observe that
$$
(f_n^{m+n})_*\big( (x_{m+n})^{\frac{1}{\ell^{n+m}}\al}\big)=(x_n)^{\frac{1}{\ell ^n}\al}.
$$
\medskip

Notice that $(x_n)^{\frac{1}{\ell ^n}(-\al)}\neq ((x_n)^{{\frac{1}{\ell ^n}}\al})^{-1}.$

\bigskip

\noindent
{\bf Lemma 2.0.} Let $z$ be a $\Qbb$-point or a tangential point defined over $\Qbb$ of $\Pbb ^1 \setminus \{0,1,\infty \}$.
\begin{enumerate}
 \item[A)] Let $\ga$ be a path on 
$\Pbb ^1 _{\bar \Qbb } \setminus \{0,1,\infty \}$ from $\01$ to $z$. Then there is a compatible family of paths 
\[
 (\ga _n)_{n\in\Nbb }\in {\varprojlim}\, \pi (V_n;\ga _n(1) , \01) 
\]
such that
\begin{enumerate}
 \item [i)] 
\[
 \ga _0=\ga\; ;
\]
 \item [ii)] if $z$ is a $\Qbb$-point then $\big(\ga _n(1)\big)_{n\in\Nbb }$ is a compatible family of $\ell^n$-th roots of $z$;
 \item [iii)] if $z$ is a tangential point then $\big(\ga _n(1)\big)_{n\in\Nbb }$ is a compatible family of  tangential points, i.e. 
$f_n^{m+n}\big( \ga _{n+m}(1)\big) =\ga _n(1)$ for all $n$ and $m$;
 \item [iv)] the compatible family of paths $(\ga _n)_{n\in \Nbb} $ is uniquely determined by the path $\ga$.
\end{enumerate}
 \item[B)]  Let us assume that a compatible family $(z^{\frac{1}{\ell  ^n}})_{n\in\Nbb}$ of $\ell^n$-th roots of $z$ is given or that a compatible family of 
tangential points is given. Then there exists a compatible family of paths
 \[
  (\ga _n)_{n\in\Nbb }\in {\varprojlim}\, \pi (V_n; z^{\frac{1}{\ell  ^n}} , \01) \;.
 \]
  \item [C)] Let $(z^{\frac{1}{\ell  ^n}})_{n\in\Nbb}$ be a given compatible family of  $\ell^n$-th roots of $z$ or a given compatible family of 
tangential points lying over $z$.
Let $\ga$ be a path on 
  $\Pbb ^1 _{\bar \Qbb } \setminus \{0,1,\infty \}$ from $\01$ to $z$. Then there is $\al \in \zl $ such that a compatible family of $\ell^n$-th roots of $z$ 
  or a compatible family of tangential points lying over $z$ determined by the path 
\[
 \delta :=\ga \cdot x^\al
\]
by the homotopy lifting property for coverings is the given family $(z^{\frac{1}{ \ell  ^n}})_{n\in\Nbb}$ of  $\ell^n$-th roots of $z$ or the given compatible family of 
tangential points lying over $z$. 
\end{enumerate}
\medskip

\noindent
{\bf Proof.} If $\ga$ is a path on $\Pbb ^1(\Cbb ) \setminus \{0,1,\infty \}$ then the existence and the uniqueness of the compatible family $(\ga _n)_{n\in \Nbb} $ 
follows from the uniqueness of the homotopy lifting property for coverings. 
If $\ga$ is arbitrary then we use the fact that the set $\pi (\Pbb ^1(\Cbb ) \setminus \{0,1,\infty \};z,\01 )$ is dense in  
$\pi (\Pbb ^1 _{\bar \Qbb } \setminus \{0,1,\infty \};z,\01 )$. 
The points ii) and iii) of A are clear.

To show the point B) of the lemma observe that the profinite sets $\pi (V_n;z^{\frac{1}{\ell  ^n}},\01 )$ are compact and the maps 
\[
 (f_n^{n+1})_* :\pi (V_{n+1};z^{\frac{1}{\ell  ^{n+1}}},\01 )\to \pi (V_n;z^{\frac{1}{\ell  ^n}},\01 )
\]
are continuous. Therefore the set ${\varprojlim}\, \pi (V_n; z^{\frac{1}{ \ell  ^n}} , \01)$ is not empty. 
Hence we get a compatible family of paths. In fact we get infinite many of compatible families.

It rests to show C). Lifting the path $\ga$ to the coverings $V_n$ of $V_0$ we get a new compatible family of $\ell^n$-th roots of $z$, which we can write in the form
\[
 (\xi ^{-\al}_{\ell ^n}  z^{\frac{1}{\ell ^n}})_{n\in \Nbb}
\]
for some $\al \in \zl$. Then lifting the path $\delta :=\ga \cdot x^\al$ to the covering $V_n$ we get the given family 
$ (z^{\frac{1}{\ell ^n}})_{n\in \Nbb}$.
\hpb

\medskip

Let $z$ be a $\Qbb$-point or a tangential point defined over $\Qbb$ of $\Pbb ^1\setminus \{0,1,\infty \}$.
Let $\ga$ be a path from $\01$ to $z$. Let
\[
 (\ga _n)_{n\in\Nbb} \in \varprojlim \pi (V_n;\zn , \01)
\]
be a compatible family of paths such that $\ga _0=\ga$.

We take the Kummer character $\kappa (z)$ equal $l(z)_{\ga _0}$.
For  $\si \in G_\Qbb$, the Kummer character evaluated at $\si$, $\kappa (z)(\si)\in\Zbb _\ell$. Let us set 
$$
\ga _{n,\si}:=\big( g^{(n)}_{\kappa (z)(\si )}(\ga _n)\big) \cdot (x_n)^{{\frac{1}{\ell ^n}}\ka (z)(\si )}.
$$
Then $\ga _{n,\si}$ is a path from $\01$ to $\xi _{\ell ^n}^{\ka (z)(\si)}\zn$. For each $n$ we have
$$
(f_n^{n+1})_*(\ga _{n+1,\si})=\ga _{n,\si}.
$$
Hence it follows that
$$
(\ga _{n,\si})_{n\in \Nbb }\in {\varprojlim}\, \pi (V_n;\xi _{\ell ^n}^{\ka (z)(\si)}\zn , \01).
$$

\medskip

\noindent
{\bf Definition 2.1.} Let $\si \in G_\Qbb $. Let us set
$$
\dfk _{\ga _n}(\si):=\ga _{n,\si}^{-1}\cdot \si (\ga _n)\in \pi _1(V_n,\01 ) 
$$
and
$$
\De _{\ga _n}(\si ):=E_n(\ga _{n,\si}^{-1}\cdot \si (\ga _n))\in \Qbb \{\{\Ybb _n\}\}.
$$

\medskip

For $n=0$ we get
$$
\De _{\ga _0}(\si)=\exp (-\ka (z)(\si)X_0)\cdot E_0(\ga_0^{-1}\cdot \si (\ga _0))=\exp (-\ka (z)(\si )X_0)\cdot \Lambda _{\ga _0}(\si ).
$$

Observe that
\begin{equation} \label{eq:compDelta}
(f_n^{m+n})_*(\De _{\ga _{m+n}}(\si ))=\De _{\ga _n}(\si ).
\end{equation}
 

\bigskip

We denote by \[\Mc _n\]  the set of all monomials in non-commuting variables belonging to $\Ybb _n$.

\medskip

\noindent{\bf Definition 2.2.} Let $z$ be a $\Qbb$-point of $\Pbb ^1\setminus \{ 0,1,\infty \}$ or a tangential point defined over $\Qbb $. 
Let $\ga$ be a path from $\01$ to $z$ on $V_0$.
Let 
$
(\ga _n)_{n\in\Nbb }\in {\varprojlim}\, \pi (V_n;\zn , \01) 
$
be such that $\ga _0=\ga$.
The functions $$\lambda _w^n(z)\;\;\;{\rm  and}\;\;\;li_w^n(z)$$ on $G_\Qbb$ are defined by the following equalities
$$
\De _{\ga _n}(\si ) =1+\sum _{w\in \Mc _n}\lambda _w^n(z)(\si )\cdot w
$$
and
$$
\log \De _{\ga _n}(\si ) =\sum _{w\in \Mc _n}li _w^n(z)(\si )\cdot w.
$$
\medskip

\noindent
For integers $0\leq i_1,i_2,\ldots ,i_r<\ell ^n$ we set \[ w(i_1,i_2,\ldots ,i_r)=Y_{n,i_1} Y_{n,i_2}  \ldots Y_{n,i_r}.\]

\bigskip

\noindent{\bf Proposition 2.3.}  Let $r>0$. The functions
$$
K_r^{(n)}(z)(\si ):(\Zbb /\ell^n)^r\to \Qbb _\ell
$$
\[
( {\rm resp.}\;\; G_r^{(n)}(z)(\si ):(\Zbb /\ell^n)^r\to \Qbb _\ell\; )
\]
defined by the formula
$$
K_r^{(n)}(z)(\si )(i_1,i_2,\ldots ,i_r):=li ^n_{w(i_1,i_2,\ldots ,i_r)}(z)(\si ) 
$$
\[
( {\rm resp.}\;\;G_r^{(n)}(z)(\si )(i_1,i_2,\ldots ,i_r):=\lambda ^n_{w(i_1,i_2,\ldots ,i_r)}(z)(\si ) \; ),
\]
where $0\leq i_1,i_2,\ldots ,i_r<\ell ^n$ 
define  a measure 
\[K_r(z)(\si )=\big( K_r^{(n)}(z)(\si )\big) _{n\in \Nbb }\]
\[
 ( {\rm resp.}\;\;G_r(z)(\si )=\big( G_r^{(n)}(z)(\si )\big) _{n\in \Nbb }\; )
\]
 on $(\Zbb _\ell)^r$ with values in $\Qbb _\ell$.

\medskip

\noindent {\bf Proof.}  It follows from the formulae \eqref{eq:fonYn} and \eqref{eq:compDelta}  that $K_r(z)(\si )$ and $G_r(z)(\si )$  are distributions on 
$(\Zbb _\ell)^r$.  Both  distributions are bounded because we are in the fixed degree $r$ and therefore the denominators cannot be worse than $(r!)^r$. \hpb 

\medskip

We denote by \[    d_r \] the smallest positive integer such that the measures $K_r(z)(\si)$ and  $G_r(z)(\si)$ have values in $\ell ^{-d_r}\Zbb _\ell $.

\medskip

\noindent
Below we point out some elementary  properties of the measures $K_r(z)(\si)$. To simplify the notation we shall omit $\si$ and write $K_r(z),\; l(z)$, $li_k(z),\ldots $ instead of 
$K_r(z)(\si)$, $l(z)(\si )$, $li_k(z)(\si ), \ldots $
unless it is necessary to indicate $\si$.

\medskip

\noindent
{\bf Proposition 2.4.}
\begin{enumerate}
\item[i)]  We have
$$
\int _{\Zbb _\ell}dK_1(z)=l_1 (z)_{\ga_0}\;\;{\rm and}\;\; \int _{(\Zbb _\ell)^r}dK_r(z)=0\;\;{\rm for}\;\; r>1.
$$
Let $0\leq a_1,\ldots ,a_r<\ell ^n$. Then
\[
 \int _{(a_1,\ldots ,a_r)+\ell ^n(\Zbb _\ell )^r}dK_r(z)=K_r^{(n)}(z)(a_1,\ldots ,a_r)\,.
\]

\item[ii)] The measure $\ell ^{d_r}K_r(z)\in  \Zbb _\ell[[(\Zbb _\ell)^r]]$ corresponds to the power series
$$
P(\ell ^{d_r}K_r(z))(A_1,\ldots ,A_r)=
$$
\[
\sum _{n_1=0}^\infty\ldots \sum _{n_r=0}^\infty \big( \int _{ (\Zbb _\ell)^r }C_{n_1}^{x_1}  C_{n_2}^{x_2}\ldots C_{n_r}^{x_r}d(\ell ^{d_r}K_r(z))\big)
 A_1^{n_1}  A_2^{n_2}\ldots A_r^{n_r}\, .
\] 
\item[iii)] We have 
\[
F(K_r(z))(X_1,\ldots ,X_r)=
\]
\[
\sum _{n_1=0}^\infty\ldots \sum _{n_r=0}^\infty {\frac{1}{n_1!  n_2!\ldots n_r!}}
\big( \int _{  (\Zbb _\ell)^r }  x_1^{n_1}   x_2^{n_2}\ldots  x_r^{n_r}dK_r(z)\big) X_1^{n_1} X_2^{n_2}\ldots X_r^{n_r}
\] 
in  $\Qbb _\ell[[X_1,X_2\ldots X_r]]$.
\end{enumerate}

\medskip

We recall that $z$ is a $\Qbb$-point of $\Pbb ^1\setminus \{0,1,\infty \}$ or a tangential point defined over $\Qbb$.
We recall that  $\ga :=\ga _0$ is a path on $V_0=\Pbb ^1_{\bar \Qbb }\setminus \{0,1,\infty \}$ from $\01$ to $z$. To simplify the notation 
we denote $X_0$ by $X$ and $Y_{0,0}$ by $Y$. Accordingly to Definition 2.2 we have
$$
\log \De _\ga =\sum _{w\in \Mc _0}li _w^0(z)\cdot w\;\;\;{\rm and}\;\;\;\De _\ga =1+\sum _{w\in \Mc _0}\lambda  _w^0(z)\cdot w\, .
$$
In \cite{NW} there are calculated coefficients $li _{YX^{n-1}}^0(z)$ of $\log \De _\ga$. Our next theorem  generalizes the result from \cite{NW}.

\medskip

\noindent{\bf Theorem 2.5.} Let $z$ be a $\Qbb$-point of $\Pbb ^1 \setminus \{0,1,\infty \}$ or a tangential point defined over $\Qbb $. 
Let $\ga$ be a path from $\01$ to $z$ on $\Pbb _{\bar  \Qbb}^1\setminus \{ 0,1,\infty \}$. 
Let 
$
(\ga _n)_{n\in\Nbb }\in {\varprojlim}\, \pi (V_n;\zn , \01) 
$
be a compatible family of paths such that
$\ga =\ga _0$.  
Let $$w=X^{a_0}YX^{a_1}YX^{a_2}Y\ldots X^{a_{r-1}}YX^{a_r}.$$ Then we have
\begin{equation} \label{eq:integralform}
li _w^0(z) = 
\end{equation}
\[
 \big( \prod _{i=0}^r a_i!\big) ^{-1}\int _{ (\Zbb _\ell)^r}(-x_1)^{a_0} (x_1-x_2)^{a_1}  (x_2-x_3)^{a_2}\ldots (x_{r-1}-x_r)^{a_{r-1}}  x_r^{a_r}dK_r(z) 
\] 
and
\begin{equation} \label{eq:integral formG}
\lambda _w^0(z)= 
\end{equation}
\[
 {\frac{1}{a_0! a_1! \ldots a_r!}}\int _{ (\Zbb _\ell)^r}(-x_1)^{a_0} (x_1-x_2)^{a_1}\ldots (x_{r-1}-x_r)^{a_{r-1}} x_r^{a_r}dG_r(z). 
\]

\medskip

\noindent {\bf Proof.} It follows from the formula \eqref{eq:compDelta}  that for any $n$ we have
$$
(f_0^n)_* (\log \De _{\ga _n}  )=\log \De _\ga .
$$
The term
$$
li _w^0(z)\,X^{a_0}YX^{a_1}Y\ldots X^{a_{r-1}}YX^{a_r}
$$
is one of the terms of the power series $\log \De _\ga $. We must see what terms of the power series $\log \De _{\ga _n} (\si )$, after applying $(f_0^n)_*$, contribute to the coefficient at $w$ of the power series $\log \De _\ga $. Let
\[
 w(i_1,i_2,\ldots ,i_r)=Y_{n,i_1} Y_{n,i_2}\ldots Y_{n,i_r} .
\]
It follows from \eqref{eq:fonYn} that the 
term
$$
li ^n_{w(i_1,i_2,\ldots ,i_r)}(z)Y_{n,i_1}Y_{n,i_2}\ldots Y_{n,i_r}
$$
is mapped by $(f_0^n)_*$ onto
$$
li ^n_{w(i_1,i_2\ldots i_r)}(z)\big( \exp (-i_1X)\cdot Y \cdot \exp (i_1X)\big) \cdot \big( \exp (-i_2X)\cdot Y \cdot
$$
$$
 \exp (i_2X)\big) \ldots \big( \exp (-i_rX)\cdot Y \cdot \exp (i_rX)\big) .
$$
Hence these terms  contribute  to the coefficient at $w$ of the power series $\log \De _\ga$ by the expression 
\begin{equation}\label{eq:sumformula}
\sum _{i_1=0}^{\ell^n-1}\sum _{i_2=0}^{\ell^n-1}\ldots \sum _{i_r=0}^{\ell^n-1}li^n _{w(i_1,i_2\ldots i_r)}(z) {\frac{(-i_1)^{a_0}}{a_0!}}  
{\frac{(i_1-i_2)^{a_1}}{a_1!}}\ldots {\frac{(i_{r-1}-i_r)^{a_{r-1}}}{a_{r-1}!}}  {\frac{(i_r)^{a_r}}{a_r!}}. 
\end{equation}
There are also terms with $X_n$ which contribute. But we have $(f_0^n)_*(X_n)=\ell ^nX.$ Therefore the contribution from terms containing $X_n$ tends to $0$ if $n$ tends to $\infty$.
Observe that if $n$ tends to $\infty$ then the sum \eqref{eq:sumformula} tends to the integral \eqref{eq:integralform}. \hpb

\medskip
The measures $K_r(z)$, $G_r(z)$, the functions $li_w^0(z)$, $ \lambda _w^0(z)$, $li_w^n(z)$, $ \lambda _w^n(z)$ depend on the path $\ga$,
hence we shall denote them also by $K_r(z)_\ga$,  $G_r(z)_\ga$, $li_w^0(z)_\ga$, $\lambda _w^0(z)_\ga$, $li_w^n(z)_\ga$, $ \lambda _w^n(z)_\ga$.

\bigskip

Throughout this paper we are working over $\Qbb$ though without any problems the base field $\Qbb$ can be replaced by any number field $K$. 
Only in Section 5 in the last two propositions and in Sections 10 and 11 the base field is $\Qbb (\mu _m)$.

\bigskip

\section{Inclusions}

\smallskip

In this section and in the next two sections we shall study symmetries of the measures $K_r(z)$. The symmetries considered are inclusions, rotations and 
the inversion. The symmetry relations are special cases of functional equations studied in \cite{W2}, \cite{W5} and recently in \cite{NW2} and \cite{NW3}.

\medskip

The inclusion
\[
 \iota _n^{p+n}:V_{p+n}\to V_n
\]
induces morphisms of fundamental groups 
\[
 (\iota _n^{p+n})_*:\pi _1(V_{p+n},\01 )\to \pi _1(V_n,\01 ) 
\]
and maps of torsors of paths
\[
 (\iota _n^{p+n})_*:\pi  (V_{p+n};z,\01 )\to \pi  (V_n;z,\01 )\,.
\]
The morphisms $ (\iota _n^{p+n})_*$ of fundamental groups induce morphisms of $\Qbb _\ell$-algebras
\[
 (\iota _n^{p+n})_*:\Qbb _\ell\{\{\Ybb _{p+n}\}\}\to \Qbb _\ell\{\{\Ybb _n\}\}\,.
\]
All these maps are compatible with the actions of $G_\Qbb $. Observe that
\begin{equation} \label{eq:inclusion}
 (\iota _n^{p+n})_*(X_{p+n})=X_n,\;\;(\iota _n^{p+n})_*(Y_{p+n,i})=0\;\;{\rm if}\;\; i\not\equiv 0\;\;{\rm mod}\;\; \ell^p
\end{equation}
\[
 \;\; {\rm and}\;\; (\iota _n^{p+n})_*(Y_{p+n,\ell^pi})=Y_{n,i}.
\]

Let \[ (\ga _n)_{n\in \Nbb} \in \varprojlim\, \pi (V_n;z^{1/\ell^n},\01 ) \] and for any $\si \in G_\Qbb $, let \[ (\ga _{n,\si })_{n\in \Nbb} \in \varprojlim\, \pi (V_n;\xi _{\ell ^n}^{\kappa (z)(\si )}z^{1/\ell^n},\01 ) \] be as in Section 2.

Let $M$ be a fixed natural number. It follows from the equality
\[
 f_n^{n+1}\circ \iota_{n+1}^{M+n+1}=\iota_n^{M+n}\circ f_{M+n}^{M+n+1}
\]
that the following diagram commutes

$$
\begin{matrix}
 \pi (V_{M+n+1};(z^{1/\ell^M})^{1/\ell^{n+1}},\01) &{\overset{(\iota_{n+1}^{M+n+1})_* }{\lra}} & \pi (V_{n+1};(z^{1/\ell^M})^{1/\ell^{n+1}},\01)    \\ \\
{(f_{M+n}^{M+n+1})_* }\Bigl\downarrow &&{(f_n^{n+1})_*}\Bigl\downarrow   \\ \\
 \pi (V_{M+n};(z^{1/\ell^M})^{1/\ell^{n}},\01) &{\overset{(\iota_{n+1}^{M+n})_* }{\lra}} & \pi (V_{n};(z^{1/\ell^M})^{1/\ell^{n}},\01 ) 
\end{matrix}
$$
as well as the analogous diagram of fundamental groups
$$
\begin{matrix}
 \pi _1 (V_{M+n+1} ,\01) &{\overset{(\iota_{n+1}^{M+n+1})_* }{\lra}} & \pi _1(V_{n+1} ,\01)    \\ \\
{(f_{M+n}^{M+n+1})_* }\Bigl\downarrow &&{(f_n^{n+1})_*}\Bigl\downarrow   \\ \\
 \pi _1(V_{M+n} ,\01) &{\overset{(\iota_{n+1}^{M+n})_* }{\lra}} & \pi _1 (V_{n} ,\01 ).
\end{matrix}
$$ 

Let us set 
\[
 \al _n=(\iota_n^{M+n})_*(\ga _{M+n})\;\;{\rm and}\;\;\al _{n,\si}=(\iota_n^{M+n})_*(\ga _{M+n,\si}).
\]
Observe that $\al _0$ (resp. $\al _{0,\si }$) is a path on $V_0=\Pbb _{\bar \Qbb }\setminus \{0,1,\infty \}$ from $\01$ to $z^{1/\ell^M}$ (resp. to $\xi _{\ell^M}^{\kappa (z)(\si )}z^{1/\ell^M}$).

We define 
\[
 \dfk _{\al _n}(\si )=\al _{n,\si }^{-1}\cdot \si (\al _n)\in \pi _1(V_n,\01 )
\]
and 
\[
 \Delta _{\al _n}(\si )=E_n\big( \al _{n,\si }^{-1}\cdot \si (\al _n)\big) \in \Qbb _\ell\{\{ \Ybb _n\}\}.
\]
One shows that
\[
 (f_n^{m+n})_*(\Delta _{\alpha _{m+n}}(\si ))=\Delta _{\al _n}(\si ).
\]
We define functions 
\[
 li_w^n(z^{1/\ell ^M})\;\;{\rm and}\;\;\lambda ^n_w(z^{1/\ell ^M})
\]
on $G_\Qbb$ by the equalities
\[
 \log \Delta _{\al _n}(\si )=\sum _{w\in \Mc _n}li_w^n(z^{1/\ell ^M})\cdot w\; {\rm and}\;\Delta _{\al _n}(\si )=1+\sum _{w\in \Mc _n}\lambda _w^n(z^{1/\ell ^M})\cdot w.
\]
If $z=\01$ then $z^{1/\ell ^M}$ we replace by ${\frac{1}{\ell ^M}}\10$.

Then as in Section 2 we get measures 
\[
K_r(z^{1/\ell ^M})\;\;{\rm and}\;\;G_r(z^{1/\ell ^M})\;\; {\rm on}\;\;(\Zbb _\ell)^r.
\] 
The analogue of Theorem 2.5 holds for the power series $\Delta _{\al _0}(\si )$ and  $\log \Delta _{\al _0}(\si )$.

\medskip

\noindent{\bf Theorem 3.1.}  Let $z$ be a $\Qbb$-point of $\Pbb ^1 \setminus \{0,1,\infty \}$ or a tangential point defined over $\Qbb $. 
Let $w=X^{a_0}YX^{a_1}YX^{a_2}Y\ldots X^{a_{r-1}}YX^{a_r}$. Then we have
\begin{equation} \label{eq:integral form}
li _w^0(z^{1/\ell ^M})=
\end{equation} 
\[
 {\frac{1}{a_0!  a_1! \ldots a_r!}}\int _{ (\Zbb _\ell)^r}(-x_1)^{a_0}  (x_1-x_2)^{a_1} \ldots (x_{r-1}-x_r)^{a_{r-1}}  x_r^{a_r}dK_r(z^{1/\ell ^M}) 
\]
and
\begin{equation} \label{eq:integral formG1}
\lambda _w^0(z^{1/\ell ^M})=
\end{equation}
\[
 {\frac{1}{a_0!  a_1! \ldots a_r!}}\int _{ (\Zbb _\ell)^r}(-x_1)^{a_0}  (x_1-x_2)^{a_1}\ldots (x_{r-1}-x_r)^{a_{r-1}}  x_r^{a_r}dG_r( z^{1/\ell ^M}). 
\]

\bigskip

The next result shows the relation between measures $K_r(z)$ and $K_r(z^{1/\ell ^M})$.

\smallskip

\noindent{\bf  Proposition 3.2.}  Let $z$ be a $\Qbb$-point of $\Pbb ^1 \setminus \{0,1,\infty \}$. Then we have
\[
 K_r^{(M+n)}(z)(\ell ^Mi_1,\ell ^Mi_2,\ldots ,\ell ^Mi_r)=K_r^{(n)}(z^{{\frac{1}{\ell ^M}}})(i_1,i_2,\ldots ,i_r)
\]
and
\[
 G_r^{(M+n)}(z)(\ell ^Mi_1,\ell ^Mi_2,\ldots ,\ell ^Mi_r)=G_r^{(n)}(z^{{\frac{1}{\ell ^M}}})(i_1,i_2,\ldots ,i_r).
\]

For $z=\10$ we have
\[
K_r^{(M+n)}(\10 )(\ell ^Mi_1,\ell ^Mi_2,\ldots ,\ell ^Mi_r)=K_r^{(n)}( {{\frac{1}{\ell ^M}}}\10 )(i_1,i_2,\ldots ,i_r)\,.
\]
If $0<i_1,i_2,\ldots ,i_r<\ell ^ n$ then
\[
 K_r^{(n)}( {{\frac{1}{\ell ^M}}}\10 )(i_1,i_2,\ldots ,i_r)=K_r^{(n)}(\10 )(i_1,i_2,\ldots ,i_r)\,.
\]

\medskip

\noindent {\bf Proof.} From the very definition of paths $\alpha _n$ and $\al _{n,\si }$ we get that for each $n$
\[
 (\iota _n^{M+n})_*(\log \Delta _{\ga _{M+n}})=\log \Delta _{\al _n} .
\]
Comparing coefficients on both sides of the equality and using the equalities \eqref{eq:inclusion} we get the first two equalities of the proposition
as well as the first equality involving $\10$. The last equality follows from the fact that the path from 
$\10$ to ${\frac{1}{\ell ^M}}\10$ is in an infinitesimal neighbourhood of $1$. \hpb

\bigskip

\section{Inversion}

\smallskip

We start with the special case of the measure $K_1(\10 )$.
Let $p_n$ be the standard path from $\01$ to ${\frac{1}{\ell ^n}}\10$ on $V_n$. Let 
\[
 h:V_n\to V_n
\]
be defined by 
\[
 h(\zfk )=1/\zfk .
\]
Let $q_n:=h(p_n)^{-1}$, let $s_n$ be a path from  ${\frac{1}{\ell ^n}}\10$ to  ${\frac{1}{\ell ^n}}\2 $ as on the picture and let $\Gamma _n:=q_n\cdot s_n \cdot p_n$.

\[
 \;
\]

\[
\;
\]

\bigskip

$${\rm Picture \;4}$$

\medskip

For $\si \in G_\Qbb$, let us define coefficients $a_i(\si )$ by the congruence
\[
 \log \Lambda _{p_n}(\si )\equiv \sum _{i=0}^{\ell^n-1}a_i(\si )Y_{n,i}\;\;{\rm mod}\;\; \Gamma ^2 L(\Ybb _n).
\]
It follows from \cite{W1} that
\[
 \ffk _{\Ga _n}=\big( p_n ^{-1}\cdot s_n^{-1}\cdot q_n^{-1}\cdot (h_*\ffk _{p_n})^{-1}\cdot q_n\cdot s_n\cdot p_n \big) \cdot (p_n^{-1}\cdot \ffk _{s_n}\cdot p_n)\cdot \ffk _{p_n}\, . 
\]
Hence we get
\[
 \log \Lambda _{\Ga _n}=-\log (h_*\Lambda _{p_n})+\log \Lambda _{s_n} +\log \Lambda _{p_n}\;\;{\rm mod}\;\;\Gamma ^2 L(\Ybb _n) \, .
\]
Observe that \[ \log \Lambda _{s_n} ={\frac{\chi -1}{2}}Y_{n,0}\] 
and
\[ -\log h_*\Lambda _{p_n}\equiv -\sum _{i=1}^{\ell ^n-1}a_iY_{\ell ^n-i}\;\;{\rm mod}\;\;\Gamma ^2 L(\Ybb _n) \, .\] 
Hence it follows that 
\begin{equation}\label{eq:ai-a-i}
 \log \Lambda _{\Ga _n}\equiv {{\frac{\chi -1}{2}}}Y_{n,0}+\sum _{i=1}^{\ell^n-1}(a_i -a_{{\ell ^n}-i})Y_{n,i}\;\;{\rm mod }\;\; \Ga ^2 L(\Ybb _n).
\end{equation}
We recall that 
for $\al \in \Qbb _\ell$ and $k\in \Nbb$ we denote by $C_k^\al$ the binomial coefficients.

\medskip

\noindent{\bf Lemma 4.1.}  For $0<i<\ell ^n$ we have
\[
 a_i(\si )-a_{\ell ^n-i}(\si )=\big( {\frac{i}{\ell ^n}}-{\frac{1}{2}}\big) -
\big( \chi(\si ) {\frac{\langle i\chi(\si ) ^{-1}\rangle }{\ell ^n}}-\chi (\si ) {\frac{1}{2}}\big)=E^{(n)}_{1,\chi (\si )}(i)\, .
\]

\medskip

\noindent {\bf Proof.} Let $\zfk$ be the standard local parameter at $0$ corresponding to $\01$. Then $u=1/\zfk$ is the local parameter at $\infty$ correspoding to $\overset{\to}{\infty 1}$.
Notice that 
\[
 \ffk _{\Ga _n}\equiv \prod _{i=0}^{\ell ^n-1} y_{n,i}^{c_i}\;\;{\rm mod}\;\; \Ga ^2 \pi _1 (V_n,\01). 
\]
To calculate the coefficients $c_i$ we shall act on 
\[
 (1-\xi _{\ell ^n}^{-i}\zfk )^{1/\ell ^m}=\sum _{k=0}^{\infty}C_k^{1/\ell ^ n}(-\xi _{\ell ^n}^{-i}\zfk )^k
\]
by the path $\ffk _{\Ga _n}(\si )=\Ga _n^{-1}\cdot \si \cdot \Ga _n \cdot \si ^{-1}$. We have
\[
(1-\xi _{\ell ^n}^{-i}\zfk )^{{\frac{1}{\ell ^m}}}{\overset{\si ^{-1}}{\lra}}(1-\xi _{\ell ^n}^{-i\chi(\si ^{-1})}\zfk )^{{\frac{1}{\ell^m}}}{\overset{\Gamma _n}{\lra}}
\]

\[
({\frac{1}{\zfk }})^{-1/\ell^m}  ({\frac{1}{\zfk }}-\xi _{\ell ^n}^{-i\chi(\si ^{-1})} )^{{\frac{1}{\ell^m}}}=u^{-1/\ell^m}  
(-\xi _{\ell ^n}^{-i\chi(\si ^{-1})} )^{{\frac{1}{\ell^m}}}  (1-\xi _{\ell ^n}^{i\chi(\si ^{-1})}u)^{{\frac{1}{\ell^m}}}{\overset{\si  }{\lra}}
\]

\[
 \si \big( (-\xi _{\ell ^n}^{-i\chi (\si ^{-1})})^{1/\ell ^m}\big)   u^{-1/\ell^m}  (1-\xi _{\ell ^n}^{i}u)^{1/\ell ^m}{\overset{\Gamma _n^{-1}}{\lra}}
\]

\[
\si \big( (-\xi _{\ell ^n}^{-i\chi (\si ^{-1})})^{1/\ell ^m}\big)   (-\xi ^i_{\ell ^n})^{1/\ell ^m}  (1-\xi _{\ell ^n}^{-i}\zfk )^{1/\ell ^m}.
\]
To fix the value of
\begin{equation}\label{eq:here}
 \si \big( (-\xi _{\ell ^n}^{-i\chi (\si ^{-1})})^{1/\ell ^m}\big)   (-\xi ^i_{\ell ^n})^{1/\ell ^m}
\end{equation}
we need to prolongate by analytic continuation $(1-\xi _{\ell ^n}^{-i}z)^{1/\ell ^m}$along $\Gamma _n$ and compare with 
$u^{-1/\ell^m}  (1-\xi _{\ell ^n}^{i}u)^{1/\ell ^m}$.

We parametrize (a part of) the path $s_n$ by
\[
 [0,\pi ]\ni \phi \longmapsto 1+\epsilon e^{{\sqrt {-1}}(\pi +\phi )}\,.
\]
We get that $(1-(1+\epsilon e^{{\sqrt {-1}}(\pi +\phi )}))^{1/\ell ^m}$ tends to 
$e^{\frac {{\sqrt {-1}}\pi}{\ell ^m}}\big( {\frac{1}{1+\epsilon }}\big)^{-1/\ell ^m} \big(1-{\frac{1}{1+\epsilon }}\big)^{1/\ell ^ m}$
if $\phi $ tends to $\pi$. Therefore
\[
\big( e^{\frac {{\sqrt {-1}}\pi}{\ell ^m}}\big) ^ {-1}(1-z)^{1/\ell ^m}=u^{-1/\ell ^m}(1-u)^{1/\ell ^m}\,.
\]
Hence it follows that  
\[
 c_i=\big( {\frac{i}{\ell ^n}}-{\frac{1}{2}}\big) -
\big( \chi(\si ) {\frac{\langle i\chi(\si ) ^{-1}\rangle }{\ell ^n}}-\chi (\si ) {\frac{1}{2}}\big)\,.
\]
\hpb

\medskip

Because of the importance of the lemma we gave a second proof.

\noindent {\bf Second proof.}  Let $0\leq i<\ell $. Let
\[
 \Phi _i:(V_n,\01 )\to (V_0,{{\overrightarrow{0\xi _{\ell ^n}^{-i}}}})
\]
be given by
\[
 \Phi _i(z)=\xi _{\ell ^n}^{-i}  z\,.
\]
Then we have
\begin{equation}\label{eq:1*}
(\Phi _i)_*(p_n^{-1}\cdot \si (p_n))\equiv (\Phi _i)_*(y_i)^{a_i(\si )}\;\;{\rm mod}\;\;\Gamma ^2\pi _1(V_0, {{\overrightarrow{0\xi _{\ell ^n}^{-i}}}}).
\end{equation}
Let $t_i\in \pi  (V_0; {{\overrightarrow{0\xi _{\ell ^n}^{-i}}}},\01 )$ be as on the picture.
\[
\; 
\]

\[
 \;
\]

$$\;$$

\[
 {\rm Picture\; 5}
\]
Observe that 
\begin{equation}\label{eq:2*}
t_i^{-1}\cdot (\Phi _i)_*(x)\cdot t_i=x,\;\;\; t_i^{-1}\cdot (\Phi _i)_*(y_i)\cdot t_i=y
\end{equation}
in $\pi _1(V_0,\01 )$. For any $\si \in G_\Qbb $ and any $0\leq i<\ell ^n$ we have
\[
 (\Phi _i)_*\circ \si =\si \circ (\Phi _{\langle i\chi (\si ^{-1})\rangle})_*\, .
\]
Hence we get
\[
(\Phi _i)_*(p_n^{-1}\cdot \si (p_n))=(\Phi _i)_*(p_n^{-1})\cdot \si (\Phi _{\langle i\chi(\si ^{-1})\rangle }(p_n))\, . 
\]
To simplify the notation let us set
\[
 q_i:=(\Phi _i)_*(p_n)\:/;\;{\rm end}\;\;\; Q_i:=q_i\cdot t_i\, .
\]
Then it follows from \eqref{eq:1*} and \eqref{eq:2*} that 
\begin{equation}\label{eq:3*}
Q_i^{-1}\cdot \si (Q_{\langle i\chi (\si ^{-1})\rangle })=t_i^{-1}\cdot q_i^{-1}\cdot \si(q_{\langle i\chi (\si ^{-1})\rangle })\cdot t_i\cdot t_i^{-1}\si (t_{\langle i\chi (\si ^{-1})\rangle})\equiv y^{a_i(\si )}\cdot x^{r_i(\si )}\;\ 
\end{equation}
\[
  {\rm modulo}\;\;\Gamma ^2\pi _1(V,\01 )
\]
for some $r_i(\si )\in \Zbb _\ell$.

Let $\zfk$ be the standard local parameter at $0$ corresponding to $\01$. Then $\tfk=\xi _{\ell ^n}^{i\chi (\si ^{-1})}\zfk$ is a local parameter at $0$ 
corresponding to 
${\overrightarrow{0\xi_{\ell ^n}^{-i\chi (\si ^{-1})} }}$ and $\tfk _1=\xi _{\ell ^n}^{i}\zfk$ is a local parameter at $0$ corresponding to 
${\overrightarrow{0\xi_{\ell ^n}^{-i} }}$.
 We calculate the action of $t_i^{-1}\cdot \si \cdot t_{\langle i\chi (\si ^{-1})\rangle}\cdot \si ^{-1}$ on $\zfk ^{1/\ell ^m}$. 
We have

\medskip

\[
\zfk ^{1/\ell ^m}\; {\overset{\si ^{-1}}{\lra}}\; \zfk ^{1/\ell ^m}\; {\overset{t_{\langle i\chi (\si ^{-1})\rangle }}{\lra}}\; 
\big(\xi _{\ell ^n}^{i\chi (\si ^{-1})} \big) ^{-1/\ell ^m} \tfk ^{1/\ell ^m}
\]
\[
  {\overset{\si }{\lra}}\; \si \big(   \big(\xi _{\ell ^n}^{i\chi (\si ^{-1})} \big) ^{-1/\ell ^m}   \big)\tfk _1^{1/\ell ^m}\; {\overset{t_i ^{-1}}{\lra}}\;
\si \big(   \big(\xi _{\ell ^n}^{i\chi (\si ^{-1})} \big) ^{-1/\ell ^m}   \big)(\xi ^i _{\ell ^n})^{1/\ell ^m}\zfk ^{1/\ell ^m}\,.
\]

\medskip
Observe that $\tfk _1^{1/\ell ^m}$ (resp. $\tfk  ^{1/\ell ^m}$) is real positive on $\varepsilon \cdot  \xi _{\ell ^n}^{-i}$ 
(resp. on $\varepsilon \cdot  \xi _{\ell ^n}^{-i\chi (\si ^{-1})}$) for $\varepsilon >0$. This fixes values   
$(\xi ^i _{\ell ^n})^{1/\ell ^m}$ (resp. $(\xi _{\ell ^n}^{i\chi (\si ^{-1})} \big) ^{-1/\ell ^m}$) for $0<i<\ell ^n$ which are $\xi _{\ell ^{n+m}}^i$ (resp. 
$\xi _{\ell ^{n+m}}^{- i\chi (\si ^{-1})  }$).

Hence we get that

\[
 r_i(\si )=  {\frac{i}{\ell ^n}}  -
  \chi(\si ) {\frac{\langle i\chi(\si ) ^{-1}\rangle }{\ell ^n}}\,.
\]
Let $h:V_0\to V_0$ be defined by
\[
 h(z)=1/z.
\]
The path $\Gamma _0$ on $V_0$ we denote by $\Gamma$.
Observe that
\begin{equation}\label{eq:4*}
 \Gamma ^{-1}\cdot h_*(x)\cdot \Gamma =y^{-1}\cdot x^{-1},\;\;\;\Gamma ^{-1}\cdot h_*(y)\cdot \Gamma =y
\end{equation}
and
\begin{equation}\label{eq:5*}
 (h(Q_i)\cdot \Gamma )^{-1}\cdot Q_{-i}=y^{-1}\cdot x^{-1}\, .
\end{equation}
It follows from \eqref{eq:3*} and \eqref{eq:4*} that 
\[
\Ga ^{-1}\cdot h(Q_i)^{-1}\cdot h(\si (Q_{\langle i\chi (\si ^{-1})\rangle }))\cdot \Ga \equiv y^{a_i(\si )}\cdot (y^{-1}\cdot x^{-1})^{r_i(\si )}\;\;{\rm mod}\;\;\Ga ^2\pi _1(V_0,\01 ) .
\]

\medskip

On the other side it follows from \eqref{eq:5*} and \eqref{eq:3*} that

\medskip
\[
 \Ga ^{-1}\cdot h(Q_i)^{-1}\cdot h(\si (Q_{\langle i\chi (\si ^{-1})\rangle }))\cdot \Ga =
\]
\[
(h(Q_i)\cdot \Ga )^{-1}\cdot Q_{-i}\cdot (Q_{-i})^{-1}\cdot \si (Q_{\langle -i\chi (\si ^{-1})\rangle })\cdot \si (x)\cdot \si (y)\cdot (\Ga ^{-1}\cdot \si (\Ga ))^{-1}\equiv 
\]
\[
 y^{-1}\cdot x^{-1}\cdot y^{a_{-i}(\si )}\cdot x^{r_{-i}(\si )}\cdot x^{\chi (\si )}\cdot y^{\chi (\si )}\cdot y^{-{\frac{1}{2}}(\chi (\si )-1)}\;\;{\rm mod}\;\;\Ga ^2\pi _1(V_0,\01 ) .
\]
Hence comparing the right hand sides of both congruences we get
\[
 a_i(\si )-r_i(\si )=a_{-i}(\si )+{\frac{1}{2}}(\chi (\si )-1)\, .
\]
Therefore we have
\[
 a_i(\si )-a_{-i}(\si )=r_i(\si )+{\frac{1}{2}}(\chi (\si )-1)=E_{1,\chi (\si )}(i)\, .
\]
\hpb

In \cite{NW2} there is still another proof of Lemma 4.1.
In the second part of the paper  we shall consider general case.


\medskip

\section{Measures $K_1(z)$}

\smallskip

In this section we present some elementary properties of measures $K_1(z)$. Most of these properties are already well known and we just collect them.

If $\mu$ is a measure on $\Zbb _\ell$ we denote by $\mu ^\times $ the restriction of $\mu$ to $\Zbb _\ell ^\times$, i.e. 
\[
 \mu ^\times =i^! \mu,
\]
where $i:\Zbb _\ell ^\times \hookrightarrow \Zbb _\ell$ is the inclusion. 

We define
\[
 m(n):\Zbb _\ell  \to \Zbb _\ell
\]
by the formula $m(n)(x)=\ell ^n x$.

\medskip

\noindent{\bf Proposition 5.1.} Let $z$ be a $\Qbb$-point of $\Pbb ^1\setminus \{ 0,1,\infty \}$. 
Let $\ga$ be a path from $\01$ to $z$. The measure $K_1(z)$ associated with the path $\ga$ from $\01$ to $z$ has the following properties:
\begin{enumerate}
 \item[i)]  
 \[ 
  F(K_1(z))(X)=\sum _{k=0}^\infty li_{k+1}(z) _\ga \cdot X^{k}\,;
  \]
 \item[ii)]  
  \[ 
  P(K_1(z))(A)=\sum _{k=0}^\infty t_{k+1}(z) _\ga \cdot A^{k}\,;
  \]
 \item[iii)] 
 \[ 
 m(n)^!K_1(z)=K_1(z^{\frac{1}{\ell ^n}})\,;  
 \]
 \item[iv)]  
 \[
 \int _{\ell ^n\Zbb _\ell}dK_1(z)=l(1-z^{1/\ell ^n })_{\alpha _0}\,;
 \]
 \item[v)]   
\[
 \int _{\Zbb _\ell}x^mdK_1(z)=\sum _{k=0}^\infty \ell ^{km}\int _{\Zbb _\ell ^\times }x^mdK_1(z^{1/\ell ^k})^\times \;\;{\rm for}\;\; m\geq 1\; .
\]

\end{enumerate}

\medskip

\noindent
{\bf Proof.} It follows from \eqref{eq:F(mu)} that 
\[
 F(K_1(z))(X)=\sum _{k=0}^\infty {\frac{1}{k!}}\big( \int _{\Zbb _\ell}x^ndK_1(z)\big) X^k\,.
\]
Observe that
\[
 li_{k+1}(z)_\ga=li^{0}_{YX^{k }}(z)_\ga \;.
\]
Hence it follows from Theorem 2.5 that
\[
 li_{k+1}(z)_\ga={\frac{1}{k!}}\int_{\Zbb_\ell}x^kdK_1(z)\;\;\;{\rm for}\;\;\;k\geq 0\;.
\]
Therefore we get the formula i) of the proposition.

We recall that the functions $t_n(z)_\ga$ are defined by the congruences \eqref{eq:congruencefort}. 
We embed   the group  $\pi _1(\Pbb ^1_{\bar \Qbb }\setminus \{0,1,\infty \},\01 )$ into $\Zbb _\ell \{\{A,B\}\}^\times$ sending $x$ to $1+A$ and $y$ to $1+B$. 
Then the image of $x^{-l(z)_\ga }\cdot \ffk  _\ga $ is the formal power series
\[
 1+\sum _{k=0}^\infty t_{k+1}(z)_\ga B\cdot A^{k}+\ldots \, ,
\]
where we have written only terms with exactly one $B$ and which start with $B$. 
Substituting $\exp X $ for $1+A$ and $\exp Y $ for $1+B$ we get the formal power series
\begin{equation}\label{eq:745}
 (\exp (-l(z)_\ga X)\cdot \Lambda _\ga (X,Y))=1+\sum _{k=0}^\infty li_{k+1}(z)_\ga YX^{k}+\ldots \,,
\end{equation}
because taking the logarithm of this power series does not change terms of degree $1$ with respect to  $Y$.
Observe that the terms on the right hand side of the formula \eqref{eq:745}, which start with $Y$ and of degree $1$ in $Y$ can be written $Y\cdot F(K_1(z))(X)$.
By the very definition we have
\[
  F(K_1(z))(X)=P(K_1(z))(\exp X-1)\;.
\]
Hence it follows that
\[
  P(K_1(z))(A)=\sum _{k=0}^\infty t_{k+1}(z) _\ga \cdot A^{k}\,.
\]

\smallskip

Let $0\leq i<\ell ^ M$. Then we have 
$K_1(z^ {1/\ell ^ n}) (i+\ell ^ M\Zbb _\ell )=K_1^ {(M)}(z^ {1/\ell ^ n})(i)=K_1^ {(M+n)}(z)(\ell ^ ni)$ by Proposition 3.2. Calculating farther we get
$K_1^ {(M+n)}(z)(\ell ^ ni)=K_1(z)(\ell ^ ni+\ell ^ {M+n}\Zbb _\ell )=K_1(z)(m(n)(i+\ell ^ M\Zbb _\ell))   
=$ $m(n)^ !K_1(z)( i+\ell ^ M\Zbb _\ell).$ Hence we have shown the point iii).

To show the point iv) observe that Proposition 3.2 implies 
\[
 \int _{\ell ^ n\Zbb _\ell}dK_1(z)=K_1^ {(n)}(z)(0)=K_1^ {(0)}(z^ {1/\ell ^ n})(0) \; . 
\]
Notice that $K_1^ {(0)}(z^ {1/\ell ^ n})(0)  $ is the coefficient at $Y$ of the element $\Delta _{\alpha _0}$, 
hence it is equal $li _Y^ {(0)}(z^ {1/\ell ^ n})  =l(1-z^ {1/\ell ^ n})_{\alpha _0} $. (We recall that $\al _0$ is $\ga _n$ 
considered on $\Pbb ^1_{\bar \Qbb}\setminus \{0,1,\infty \}$.)

To prove the point v) we present $\Zbb _\ell$ as the following finite disjoint union of compact-open subsets 
\[
 \zl =\zlt \cup\ell \zlt\cup \ldots \cup \ell ^{n-1}\zlt \cup \ell ^n\zl \;.
\]
Observe that 
\[
 \int _{\ell ^k\zlt }x^mdK_1(z)=\int _{\zlt}(\ell ^kx)^md(m(k)!K_1(z))
\]
by the formula \eqref{eq:fcirc phi2}. It follows from the point iii) already proved that
\[
 \int _{\zlt}(\ell ^kx)^md(m(k)^!K_1(z))=\ell ^{km}\int _{\zlt}x^mdK_1(z^{1/\ell ^k})\; .
\]
Hence we get that
\[
  \int _{ \zl }x^mdK_1(z)=\sum _{k=0}^{n-1}\ell ^{km}\int _{\zlt}x^m dK_1(z^{1/\ell ^k}) +\ell ^{nm}\int _{\zl }x^m dK_1(z^{1/\ell^n})\; .
\]
Observe that the term $\ell ^{nm}\int _{\zl }x^m dK_1(z^{1/\ell^n})$ tends to $0$ if $n$ tends to $\infty$. Hence we have
\[
  \int _{ \zl }x^m dK_1(z)=\sum _{k=0}^{\infty}\ell ^{km}\int _{\zlt}x^m dK_1(z^{1/\ell ^k})\; .
\]
\hpb

\medskip

In the next proposition we indicate properties of the measure $K_1(\10 )$.

\smallskip

\noindent{\bf Proposition 5.2.}  Let $p$ be the standard path on $\Pbb ^1_{\bar \Qbb}\setminus \{0,1,\infty\}$ from $\01$ to $\10$. Let $K_1(\10 )$ be the measure associated with the path $p$. We have 
\begin{enumerate}
 \item[i)] \[   \big( m(n)^!K_1(\10 )\big) ^\times =K_1(\10 )^\times\;;\]
 \item[ii)] 
 \[
 \int _{\Zbb _\ell}dK_1(\10 )=0\;\;{\rm and}\;\; \int _{\ell ^n\Zbb _\ell}dK_1(\10 )=\kappa ({\frac{1}{\ell ^ n}})  \;\;{\rm for}\;\; n>0\; ;
 \]
 \item[iii)] 
 \[
  \int _{\zl}x^k dK_1(\10 )={\frac{1}{1-\ell ^k}} \int _{\zlt}x^k dK_1(\10 )\;.
 \]
\end{enumerate}
 
\medskip

\noindent
{\bf Proof.} The lifting of the path $p=p_0$ to $V_n$ is the path $p_n$ from $\01$ to ${\frac{1}{\ell ^n}}\10 $. We have
\[
 \big( m(n)^!K_1(\10 )\big) (i+\ell ^M\Zbb _\ell )=K_1(\10 )(\ell ^n i+\ell ^{M+n}\Zbb _\ell )=K_1^{(M+n)}(\10 )(\ell ^ni)\;.
\]
Observe that $K_1^{(M+n)}(\10 )(\ell ^ni)$ is the coefficient of $\log \Lambda _{p_{M+n}}$ at $Y_{M+n,\ell ^ni}$.
Assume that $\ell$ does not divide $i$. Then this coefficient is equal to the coefficient of $\log \Lambda _{p_M}$at $Y_{M,i}$, which is $K_1^{(M)}(\10 )(i)=
K_1 (\10 )(i+\ell ^M\Zbb _\ell )$. Therefore
\[
 \big( m(n)^!K_1(\10 )\big) (i+\ell ^M\Zbb _\ell ) =K_1(\10 )  (i+\ell ^M\Zbb _\ell ) 
\]
for $i$ not divisible by $\ell$. This implies the point i).

The formal power series $\Lambda _p=\Delta _p$ has no terms in degree one, hence $\int _{\zl}dK_1(\10 )=l_1(\10 )_p=0$. We have
\[
 \int _{\ell ^n\zl}dK_1(\10 )=K_1(\10 )(\ell ^n\zl )=K_1^{(n)}(0)\;.
\]
Observe that $K_1^{(n)}(0)$ is the coefficient of $\Lambda _{p_n}=\Delta _{p_n}$ at $Y_{n,0}$. 
Let $t$ be the local parametre on $V_n$ at $0$ corresponding to $\01$. 
The element $\ffk _{p_n}(\si )=p_n^{-1} \cdot \si \cdot p_n \cdot \si ^{-1}$ acts on $(1-t)^{\frac{1}{\ell ^m}}$ as follows:
\[
 (1-t)^{\frac{1}{\ell ^m}}   {\overset{\si ^{-1}}{\lra}}   (1-t)^{\frac{1}{\ell ^m}}  {\overset{p_n}{\lra}} 
({\frac{1}{\ell ^n}})^{{\frac{1}{\ell ^m}}}  s^{{\frac{1}{\ell ^m}}}  {\overset{\si  }{\lra}}
\]

\[
 \si \big( ({\frac{1}{\ell ^n}})^{{\frac{1}{\ell ^m}}}\big)   s^{{\frac{1}{\ell ^m}}} {\overset{p_n^{-1}}{\lra}} 
\si \big( ({\frac{1}{\ell ^n}})^{{\frac{1}{\ell ^m}}}\big)   \big( ({\frac{1}{\ell ^n}})^{{\frac{1}{\ell ^m}}}\big)   (1-t)^{\frac{1}{\ell ^m}}=
\xi _{\ell ^m}^{\kappa (1/\ell  ^n)}  (1-t)^{\frac{1}{\ell ^m}}\,,
\]
where $s=\ell ^n(1-t)$ is the local parametre on $V_n$ at $1$ corresponding to ${\frac{1}{\ell ^n}}\10 $.  Hence we get that
\[
 K_1^{(n)}(\10 )(0)=\kappa ({\frac{1}{\ell ^n}})
\]
and therefore $\int _{\ell ^n\Zbb _\ell}dK_1(\10 )=\kappa ({\frac{1}{\ell ^ n}}) $.

Repeating the arguments from the proof of the point v) of Proposition 5.1 we get
\[
 \int _{\Zbb _\ell }x^mdK_1(\10 )=\sum _{k=0}^\infty \ell ^{mk}\int _{\Zbb _\ell ^\times }x^mdK_1({\frac{1}{\ell ^k}}\10 )=
\sum _{k=0}^\infty \ell ^{mk}\int _{\Zbb _\ell ^\times }x^mdK_1( \10 )\;,
\]
because the measures $ K_1({\frac{1}{\ell ^k}}\10 )$ and $K_1( \10 )$ coincide on $\Zbb _\ell ^\times$. But the last series is equal 
${\frac{1}{1-\ell ^m}} \int _{\zlt}x^m dK_1(\10 )$.  \hpb   

\medskip

In the next two propositions our base field is $\Qbb (\mu _m)$.

\smallskip

\noindent
{\bf Proposition 5.3.} Let $m$ be a positive integer not divisible by $\ell$.
Let $\xi _m$ be a primitive $m$-th root of $1$. Let $\big( \xi _m^{\ell  ^{-n}}\big) _{n\in \Nbb}$ be a compatible family of $\ell ^n$-th roots of $\xi _m$ 
such that $\xi _m^{\ell  ^{-n}}\in \mu _m$ for all $n\in \Nbb$. Let $a$ be the order of $\ell$ modulo $m$. Let
\[
 (\ga _n)_{n\in\Nbb }\in {\varprojlim}\, \pi (V_n;\xi _m^{\ell  ^{-n}}   , \01) 
\]
and let $K_1(\xi _m)$ be the measure associated with the path $\ga _0$. Then we have:
\begin{enumerate}
 \item[i)]
 \[
  \int _{\Zbb _\ell}x^kdK_1(\xi _m)=\sum _{i=0}^{a-1}{\frac{\ell ^{ki}}{1-\ell ^{ka}}}\int _{\Zbb _\ell ^\times }x^kdK_1(\xi _m^{\ell ^{-i}})^\times \;\;
 {\rm for}\;\;k\geq 1\;,
 \]
 \item[ii)] 
 \[
  l_k(\xi _m^{\ell ^{-i}})_{\ga _i}=li_k(  \xi _m^{\ell ^{-i}})_{\ga _i}\;\;{\rm for}\;\;0\leq i <a\;, 
 \]
 \item[iii)] the functions
\[
 l_k(\xi _m^{\ell ^{-i}})_{\ga _i}:G_{\Qbb (\mu _m)}\to \zl (k)
\]
are cocycles for all $k$ and $0\leq i <a$.
\end{enumerate}

\medskip

\noindent
{\bf Proof.}  Observe that
\[
 \int _{\Zbb _\ell}x^ kdK_1(\xi _m)=\sum _{n=0}^ \infty \int _{\ell ^ n\Zbb _\ell ^\times}x^ kdK_1(\xi _m)=
\sum _{n=0}^ \infty \ell ^{nk} \int _{\Zbb _\ell ^ \times }x^ kdK_1(\xi _m ^ {\ell ^ {-n}})^ \times =
\]

\[
 \sum _{i=0}^ {a-1}\big( \sum _{n=0}^\infty         \ell ^ {(i+a\cdot n)k}  \int _{\Zbb _\ell ^ \times }x^ kdK_1(\xi _m ^ {\ell ^ {-i}})^ \times  \big)=\sum _{i=0}^ {a-1}
{\frac{\ell ^ {ki}}{1-\ell ^ {ka}}} \int _{\Zbb _\ell ^ \times }x^ kdK_1(\xi _m ^ {\ell ^ {-i}})^ \times\,.
\]
Hence we have shown the point i) of the proposition.

Let $0\leq i<a$. Observe that $l(\xi _m^{\ell ^{-i}})_{\ga _i}=0$ because $\ell^n$-th roots of $\xi _m^{\ell ^{-i}}$ calculated along $\ga _i$ are in $\mu _m$. 
Hence it follows that $\Lambda _{\ga _i}=\Delta _{\ga _i}$ and in consequence
\[
 l_k(\xi _m^{\ell ^{-i}})_{\ga _i}\,=\, li_k(\xi _m^{\ell ^{-i}})_{\ga _i}\;.
\]
It follows from \cite[Theorem 11.0.9]{W2} that $l_k(\xi _m^{\ell ^{-i}})_{\ga _i}$ are cocycles.   \hpb

\medskip

The last result of this section concerns distribution relations of $\ell$-adic polylogarithms.  
In \cite{NW3} we proved the following result (see also \cite[Theorem 2.1.]{W5}).

\medskip

\noindent
{\bf Theorem 5.4.} Let $m$ be a positive integer not divisible by $\ell$. 
Let $z$ be a $\Qbb$-point of $\Pbb ^1\setminus \{0,1,\infty \}$. There are $\ell$-adic paths $\ga_k$ on $\Pbb _{\bar \Qbb}^1 \setminus \{ 0,1,\infty \}$ from  
$\01$ to $\xi _m ^k$ for $k=0,1,\ldots ,m-1$ and an $\ell$-adic path 
 $\ga$ from $\01$ to $z^m$ such that
\[
 m^{n-1}  \big( \sum _{k=0} ^{m-1}li_n(\xi _m^kz)_{\ga _k}\big) \,=\,li_n(z^m)_\ga
\]
on the group $G_{\Qbb   (\mu _m)}$ for all $n\geq 1$.

\medskip

\noindent
The next result follows immediately from Theorem 2.5 and the theorem stated above.

\medskip

\noindent
{\bf Proposition 5.5.} We have the following equality of the formal power series in $\Qbb [[X]]$ 
\[
 \sum _{k=0}^{m-1}F(K_1(\xi _m^kz )_{\ga _k})(mX)=F(K_1(z^m)_{\ga})(X)\,.
\]


\section{Congruences between coefficients}

\smallskip 

Let $w=X^{a_0}YX^{a_1}Y\ldots YX^{a_r}$. In Section 2 we have shown that 

\begin{equation}\label{eq:zsec2}
 li _w^0(z)=
\end{equation}
\[
 {\frac{1}{a_0!  a_1! \ldots a_r!}}\int _{ (\Zbb _\ell)^r}(-x_1)^{a_0}  (x_1-x_2)^{a_1} \ldots (x_{r-1}-x_r)^{a_{r-1}}  x_r^{a_r}dK_r(z). 
 \]

Let $F:(\Zbb _\ell )^r \to (\Zbb _\ell )^r$ be given by $F(x_1,\ldots ,x_r)=(x_1-x_2,\ldots ,x_{r-1}-x_r,x_r)$.
Observe that $F$ is an isomorphism of $\Zbb _\ell$-modules.
It follows from the formula \eqref{eq:fcirc phi} that 
\begin{equation}\label{eq:KtobarK}
 \int _{ (\Zbb _\ell)^r}(-x_1)^{a_0}  (x_1-x_2)^{a_1}  (x_2-x_3)^{a_2}\ldots (x_{r-1}-x_r)^{a_{r-1}}  x_r^{a_r}dK_r(z)=
\end{equation}
\[
 \int _{ (\Zbb _\ell)^r}(-\sum_{i=1}^rt_i)^{a_0}  (t_1 )^{a_1}  (t_2 )^{a_2}\ldots (t_{r-1} )^{a_{r-1}}  t_r^{a_r}d(F_! K_r(z)).
\]
To simplify the notation we denote
\[
 \bar K_r(z)=F_! K_r(z).
\]

Let us decompose $(\Zbb _\ell)^r$ into a disjoint union of compact subsets
\[
 (\Zbb _\ell)^r=\bigsqcup _{n_1=0}^{\bar \infty }\ldots \bigsqcup _{n_r=0}^{\bar \infty }\big( \prod _{i=1}^r\ell ^{n_i}\Zbb _\ell ^\times\big)\, ,
\]
where bar over $\infty$ means that the summation includes $\infty$ and $\ell ^\infty \Zbb _\ell ^\times =\{ 0\}$. Observe that the subsets
\[
  \prod _{i=1}^r\ell ^{n_i}\Zbb _\ell ^\times
\]
for $n_1\neq \infty$, $n_2\neq \infty$,\ldots ,$n_r\neq \infty$ are compact-open subsets of $(\Zbb _\ell)^r$.

\medskip

Let  $n_1\neq \infty$, $n_2\neq \infty$,\ldots ,$n_r\neq \infty$. Let 
\[
 m(n_1,\ldots ,n_r):(\Zbb _\ell ^\times )^r \to (\Zbb _\ell )^r
\]
be given by
\[
 m(n_1,\ldots ,n_r)(t_1,\ldots ,t_r)=(\ell ^{n_1}t_1,\ldots ,\ell ^{n_r}t_r).
\]

\medskip

\noindent{\bf Lemma 6.1.} We have
\begin{equation}\label{eq:integralpodzielony}
\int _{\prod_{i=1}^r \ell ^ {n_i}\Zbb _\ell^\times}(-\sum_{i=1}^rt_i)^{a_0}  (t_1 )^{a_1}  (t_2 )^{a_2}\ldots (t_{r-1} )^{a_{r-1}}  (t_r)^{a_r}d\bar K_r(z)
=
\end{equation}
\[
\ell ^{\sum _{i=1}^{ r } a_in_i}\int _{ (\Zbb _\ell^\times)^r} (-\sum_{i=1}^r\ell ^{n_i}t_i)^{a_0}   (t_1 )^{a_1}  (t_2 )^{a_2}\ldots  
 (t_r)^{a_r}d\big( m(n_1,\ldots ,n_r)^!\bar K_r(z)\big) .
\]

\medskip

\noindent{\bf Proof.} The lemma  follows from the formula \eqref{eq:fcirc phi2}.    \hpb

\medskip

\noindent{\bf Lemma 6.2.} \label{lem:conv}
 Let us assume that $a_i$ are positive integers for $i=1,2,\ldots ,r$. Then we have
\begin{equation} \label{eq:conv}
\int _{ (\Zbb _\ell)^r}(-\sum_{i=1}^rt_i)^{a_0}  (t_1 )^{a_1}  (t_2 )^{a_2}\ldots (t_{r-1} )^{a_{r-1}}  (t_r)^{a_r}d\bar K_r(z)=
\end{equation}
\[
\sum _{n_1=0}^{  \infty }\ldots \sum _{n_r=0}^{  \infty } \ell ^{\sum_{i=1}^r a_in_i}\int _{ (\Zbb _\ell^\times)^r} (-\sum _{i=1}^r\ell ^{n_i}t_i)^{a_0}  
(t_1 )^{a_1}  (t_2 )^{a_2}\ldots (t_{r-1} )^{a_{r-1}}  (t_r)^{a_r}d\big( {\bf K}\big) .
\] 
where ${\bf K}=m(n_1,\ldots ,n_r)^ !\bar K_r(z)$.
\medskip

\noindent{\bf Proof.}
 Observe that for any natural number $M$ the set 
\[
 \{(n_1,n_2,\ldots ,n_r)\in \Nbb ^r\mid \sum _{i=1}^r n_ia_i<M\}
\]
is finite. This implies that the series on the right hand side of  \eqref{eq:conv}  converges. For a given $M$ 
we have the following decomposition into a finite disjoint union of compact-open subsets
\[
 (\Zb _\ell )^r=( \bigsqcup _{n_1=0}^M\ldots \bigsqcup _{n_r=0}^M\big( \prod _{i=1}^r\ell ^{n_i}\Zb _\ell ^\times \big) )\bigsqcup
 \big(  \ell ^{M+1}\Zb _\ell \big) ^r.
\]
Observe that
\[
 \int _{(\ell ^{M+1}\Zb _\ell )^r}(-\sum _{i=1}^rt_i)^{a_0}  (t_1)^{a_1}  (t_2)^{a_2}\ldots (t_r)^{a_r}d\bar K_r(z)\equiv 0\;\;{\rm mod}\;\;\ell ^{M+1-d_r}.
\]
Hence it follows from \eqref{eq:integralpodzielony} that the series on the right hand side of the equality (\ref{eq:conv}) converges to the integral on the left hand side of the equality (\ref{eq:conv}).
\hpb 

\medskip

Now we shall prove congruence relations between coefficients of the power series
\[
 \log \Delta _\ga =\sum _{w\in \Mc _0}li_w^0(z)\cdot w \in \Qbb _\ell \{\{X,Y\}\}.
\]
\medskip

\noindent{\bf Theorem 6.3. } Let $a_i$ and $b_i$ be non negative integers not divisible by $\ell$ for $i=1,2,\ldots ,r$. 
Let $w=YX^{a_1}YX^{a_2}\ldots YX^{a_r}$ and $v=YX^{b_1}YX^{b_2}\ldots YX^{b_r}$. Let $M$ be a positive integer. 
Let us assume that $a_i\equiv b_i$ modulo $(\ell -1)\ell ^M$ for $i=1,2,\ldots ,r$. Let $z$ be a $\Qbb$-point of $\Pbb ^1\setminus \{ 0,1,\infty \}$ or $z=\10$. 
Let $\ga$ be a path from $\01$ to $z$. Then for any    $\si\in G_\Qbb$ we have the following congruences between coefficients of the power series  
$ \log \Delta _\ga$ ($\log \Lambda _p$ if $z=\10$)
\[
(\prod _{i=1}^ra_i!) li _w^0(z)(\si )\equiv   (\prod _{i=1}^rb_i!) li _v^0(z)(\si ) \;\;{\rm modulo}\;\; \ell ^{M+1-d_r}\, .
\]

\medskip

\noindent{\bf Proof.} One can find $c_i\in \Zbb$ such that 
\[
 b_i=a_i+c_i(\ell -1)\ell ^M
\]
for $i=1,2,\ldots ,r$. Then for any $x\in \Zbb _\ell ^\times $ we have
\[
 x^{b_i}=x^{a_i}\cdot x^{(\ell -1)c_i\ell ^M}=x^{a_i}y ^{\ell^ M}\, ,
\]
where $y= x^{(\ell -1)c_i}\in 1+\ell \Zbb _\ell$. It implies that
\[
 x^{b_i}\equiv x^{a_i}\;\;{\rm modulo}\;\; \ell ^{M+1}
\]
for  $i=1,2,\ldots ,r$. Hence it follows that
\[
 \int _{(\Zbb _\ell ^\times )^r} t_1^{a_1}t_2^{a_2}\ldots t_r^{a_r} d(m(n_1,\ldots ,n_r)^!\bar K_r(z)(\si ))\equiv
\]
\[
 \int _{(\Zbb _\ell ^\times )^r} t_1^{b_1}t_2^{b_2}\ldots t_r^{b_r} d(m(n_1,\ldots ,n_r)^!\bar K_r(z)(\si ))\;\;{\rm modulo}\;\; \ell ^{M+1-d_r}\, .
\]
Lemma 6.2 implies that
\[
 \int _{(\Zbb _\ell   )^r }t_1^{a_1}t_2^{a_2}\ldots t_r^{a_r} d \bar K_r(z)(\si )\equiv
\]
\[
 \int _{(\Zbb _\ell )^r} t_1^{b_1}t_2^{b_2}\ldots t_r^{b_r} d \bar K_r(z)(\si )\;\;{\rm modulo}\;\; \ell ^{M+1-d_r}\, .
\] 
Therefore the theorem follows from the equality \eqref{eq:KtobarK} and Theorem 2.5.  \hpb


\medskip

\section{$\ell$-adic poly--multi--zeta functions?}

\smallskip

In this section we attempt to define non-Archimedean analogues of multi--zeta functions 
\[
 \zeta (s_1,\ldots ,s_r)=\sum _{n_1>n_2>\ldots >n_r=1}{\frac {1}{n_1^{s_1}  n_2^{s_2}\ldots n_r^{s_r}}}
\]
and poly--multi--zeta functions
\[
 \zeta _z(s_1,\ldots ,s_r)=\sum _{n_1>n_2>\ldots >n_r=1}{\frac {z^{n_1}}{n_1^{s_1}  n_2^{s_2}\ldots n_r^{s_r}}}\, .
 \]

\medskip

Let
\[
 \omega :\Zbb _\ell ^\times \to \Zbb _\ell ^\times
\]
be the Teichmuller character. If $x\in \Zbb _\ell ^\times $ we set
\[
 [x]:=x\cdot \omega (x)^{-1}\, .
\]
\medskip

\noindent{\bf Definition 7.1.} Let $0\leq \beta _i <\ell -1$ for $i=1,\ldots ,r$. Let $\bar \beta :=(\beta _1,\ldots ,\beta _r)$, 
let $\bar n:=(n_1,\ldots ,n_r)\in \Nbb ^ r$ and let 
$(s_1,\ldots ,s_r)\in (\Zbb _\ell )^r$.  
Let $z$ be a $\Qbb$-point of $\Pbb ^1 \setminus \{0,1,\infty \}$ or a tangential point defined over $\Qbb$. 
 We define
\[
 \Zc ^{\bar \beta}_{\bar n}(1-s_1,\ldots ,1-s_r;z,\si ):=
\]
\[
\int _{(\Zbb _\ell ^\times )^r}[t_1]^{s_1}t^{-1}\omega (t_1) ^{\beta _1}\ldots [t_r]^{s_1}t^{-1}\omega (t_r) ^{\beta _r}d (m (n_1,\ldots ,n_r) ^!\bar K _r(z)(\si )\,.
\]

\medskip

For $z=\10$ we should obtain $\ell$-adic non-Archimedean analogues of multi-zeta functions. However before we should divide by polynomials in $[\chi (\si )]^ s$ in order to get functions which do not depend on $\si$. 
We do not know how to do this for arbitrary $r$. Only for $r=1$ we can guess easily the required polynomial. The case $r=1$ is studied in the next section.


\medskip

\section{$\ell$-adic L-functions of Kubota-Leopoldt}

\smallskip

Now we shall consider the only case 
when we can show the expected relations of the functions constructed by us in Section 7 with the corresponding $\ell$-adic non-Archimedean functions.

We shall consider the case of $r=1$ and $z=\10$. We shall show that in this case the functions $\Zc _0^\beta (1-s;\10,\si)$ defined in 
Section 7 are in fact the  Kubota-Leopoldt L-functions multiplied by the function
\[
 s\longmapsto \omega (\chi (\si))^\beta [\chi (\si )]^s -1\, . 
\]

\medskip

We start by gathering the facts we shall need and which are crucial in identification of $\Zc _0^\beta (1-s;\10,\si)$ with the Kubota-Leopoldt L-functions.
It follows from Theorem 2.5 and the definition of $\ell$-adic Galois polylogarithms in \cite{W2} that 
\begin{equation}\label{eq:91}
 l_k(\10 )={\frac{1}{(k-1)!}}\int _{\Zbb _\ell}x^{k-1}dK_1(\10 )\, .
\end{equation}
 It follows from Proposition 5.2, point iii) that
\begin{equation}\label{eq:92}
 \int _{\Zbb _\ell}x^{k-1}dK_1(\10 )={\frac{1}{1-\ell ^{k-1}}}\int _{\Zbb _\ell^\times}x^{k-1}dK_1(\10 )\,  
\end{equation}
for $k>1$. For $k>0$ and even we have the equality
\begin{equation}\label{eq:93}
 l_k(\10 )={\frac{-B_k}{2\cdot k!}}(\chi ^k-1)
\end{equation}
(see \cite[Proposition 3.1]{W7}, another proof is in \cite{NW2}).
\medskip

In Section 7 we defined
\[
 \Zc ^\beta _0(1-s;\10 ,\si)=\int _{\Zbb _\ell ^\times }[x]^s  x^{-1} \omega (x)^\beta dK_1(\10 )(\si )\, .
\]
We shall use a modified version of the function.

\medskip

\noindent{\bf Definition 8.1.} Let $0\leq \beta <\ell -1$. Let $\si \in G_\Qbb $ be such that $\chi (\si )^{\ell -1}\neq 1$. We define 
\[
 L^\beta (1-s;\10 ,\si ):={\frac{2}{\omega (\chi (\si ))^\beta [\chi (\si )]^s-1}}\int _{\Zbb _\ell ^\times }[x]^s  x^{-1}   \omega (x)^\beta dK_1(\10 )(\si )\, .
\]

\medskip

\noindent{\bf Theorem 8.2.} Let $\si \in G_\Qbb $ be such that $\chi (\si )^{\ell -1}\neq 1$. 
\begin{enumerate}
 \item[i)] Let $k>0$ and let  $k\equiv \beta $ modulo $\ell -1$. Then we have
\begin{equation}\label{eq:11}
   L^\beta (1-k;\10 , \si )={\frac{2}{ \chi (\si )^k-1}}\int _{\Zbb _\ell ^\times }x^{k-1}  dK_1(\10 )(\si )=
{\frac{2(1-\ell ^{k-1})(k-1)!}{\chi (\si )^k-1}}l_k(\10 )(\si )\, . 
\end{equation}
 \item[ii)] Let $k>0$ and let $\beta $ be even. Then we have
\begin{equation}\label{eq:12}
  L^\beta (1-k;\10 ,\si )=-{\frac{1}{k}}B_{k,\omega ^ {\beta -k}}\,.
\end{equation}
 \item[iii)] Let $k$ and $\beta$ be even and let $k\equiv \beta $ modulo $\ell -1$. Then we have
\begin{equation}\label{eq:13}
  L^\beta (1-k;\10 ,\si )=-(1-\ell ^{k-1}){\frac{B_k}{k}}=(1-\ell ^{k-1})\zeta (1-k)\, .
\end{equation}
\end{enumerate}

\medskip

\noindent {\bf Proof.} Let us assume that $k\equiv \beta $ modulo $\ell -1$. Observe that then $[\chi (\si )]^k=\chi (\si )^k   \omega (\chi (\si ))^{-\beta}$
and $x^{k-1}=[x]^k  x^{-1}  \omega (x)^\beta$.
Hence we get
\[
   L^\beta (1-k;\10 ,\si )={\frac{2}{ \chi (\si )^k-1}}\int _{\Zbb _\ell ^\times }x^{k-1}  dK_1(\10 )(\si ).
\]
Observe that 
\[
 \int _{\Zbb _\ell ^\times }x^{k-1}  dK_1(\10 )(\si )=(1-\ell ^{k-1})\cdot (k-1)!\, l_k(\10 )(\si )\,  
\]
by the equalities \eqref{eq:92} and \eqref{eq:91}.
Now we shall prove the point ii). Let $\beta$ be even. Then we have
\[
 L^\beta (1-k;\10 ,\si ) ={\frac{2}{\omega (\chi (\si ))^{\beta -k} \chi (\si )^k-1}}\int _{\Zbb _\ell ^\times }x^{k-1}  \omega (x)^{\beta -k} dK_1(\10 )(\si )\, .
\]
It follows from Lemma 4.1 and the equality
$ E^{(n)}_{1,\chi (\si )}(\ell ^n-i)= -E^{(n)}_{1,\chi (\si )}( i)$ that 
\[
 \int _{\Zbb _\ell ^\times }x^{k-1} \omega (x)^{\beta -k} dK_1(\10 )(\si )=
{\frac{1}{2}}\int _{\Zbb _\ell ^\times }x^{k-1}  \omega (x)^{\beta -k} dE_{1,\chi(\si )}\,.
\]
Hence we get that
\[
  L^\beta (1-k;\10 ,\si )={\frac{1}{\omega (\chi (\si ))^{\beta  } [\chi (\si )]^k-1}}\int _{\Zbb _\ell ^\times }[x]^{k }  x^{-1}   \omega (x)^{\beta }
 dE_{1,\chi(\si )}\,.
\]
Therefore $ L^\beta (1-k;\10 ,\si )=-{\frac{1}{k}}B_{k,\omega ^{\beta -k}}$ by \cite[Chapter 4, Theorem 3.2.]{L}, 

It rests to show iii). If $k\equiv \beta $ modulo $\ell -1$ then
\[
  L^\beta (1-k;\10 ,\si )={\frac{2}{ \chi (\si )^k-1}}\int _{\Zbb _\ell ^\times }x^{k-1}  dK_1(\10 )(\si )
\]
by the point i) already proved. Hence it follows from \eqref{eq:91} , \eqref{eq:92} and \eqref{eq:93}  that 
\[
 {\frac{2}{ \chi (\si )^k-1}}\int _{\Zbb _\ell ^\times }x^{k-1}  dK_1(\10 )(\si )=
{\frac{2(1-\ell ^{k-1})}{ \chi (\si )^k-1}}\int _{\Zbb _\ell}x^{k-1}  dK_1(\10 )(\si )=
\]
\[
 {\frac{2(1-\ell ^{k-1})\cdot (k-1)!}{ \chi (\si )^k-1}}l_k(\10 )=-(1-\ell ^{k-1}){\frac{B_k}{k}}=(1-\ell ^{k-1})\zeta (1-k)\,.
\]
\hpb

\medskip
The $\ell$-adic $L$-functions were first defined in \cite{KL}. The other construction is given in \cite{Iw}. We shall use the definition which appear in \cite{L}. 
Following 
Lang (see \cite{L}) we define the Kubota-Leopoldt $\ell$-adic $L$-functions by
\[
 L_\ell (1-s;\Phi ):={\frac{1}{\Phi (c)[c]^s -1}}\int _{\Zbb _\ell ^\times }[x]^s \cdot x^{-1}\cdot \Phi (x)dE_{1,c}(x)\, ,
\]
where $\Phi$ is a character of finite order on $\Zbb _\ell ^\times $ and $c\in \Zbb _\ell ^\times$.

\medskip

We recall that
\begin{equation}\label{lang1}
 L_\ell (1-k,\omega ^\beta )=-{\frac{1}{k}}B_{k,\omega ^{\be -k}}
\end{equation}
for any positive integer $k$ (see \cite[Chapter 4,Theorem 3.2.]{L}). In particular if $k\equiv \be $ modulo $\ell -1$ then we have
\begin{equation}\label{lang2}
  L_\ell (1-k,\omega ^\beta )=-{\frac{1}{k}}B_{k,{\bf 1}}=-(1-\ell ^{k-1}){\frac{B_k}{k}}\;,
\end{equation}
where ${\bf 1}:\Zbb _\ell ^\times \to \{1\}$ denotes the trivial character of $\Zbb _\ell ^\times $.

\medskip

\noindent{\bf Corollary 8.3.}  Let $\beta$ be even and $0\leq \beta \leq \ell -3$.  Let $\si \in G_\Qbb $ be such that $\chi (\si ) ^{\ell -1}\neq 1$. The function 
$ L^\beta (1-s;\10 ,\si )$ does not depend on $\si $ and it is equal to the Kubota-Leopoldt $\ell$-adic $L$-function $L_\ell (1-s;\omega ^\beta )$.

\medskip

\noindent{\bf Proof.} Let $\si _1$ and $\si _2$ belonging to $G_\Qbb $ be such that  $\chi (\si _1) ^{\ell -1}\neq 1$ and $\chi (\si _2) ^{\ell -1}\neq 1$. 
It follows from the point ii) of Theorem 8.2. that 
\[
  L^\beta (1-k;\10 ,\si _1)= L^\beta (1-k;\10 ,\si _2)
\]
for $k$ a positive integer. Hence
\[
   L^\beta (1-s;\10 ,\si _1)= L^\beta (1-s;\10 ,\si _2)                                 
\]
because the functions coincide on the dense subset of $\Zbb _\ell $. 
It follows from the point iii) of Theorem 8.2 and \eqref{lang2} that $ L^\beta (1-s;\10 ,\si )$ is the Kubota-Leopoldt $\ell$-adic $L$-function  $L_\ell (1-s;\omega ^\beta )$. \hpb

\medskip

\noindent{\bf Remark 8.4.}  
\begin{enumerate}
\item[i)] If $\beta$ is odd then the functions $ L^\beta (1-s;\10 ,\si  )$ and $ \Zc ^\beta (1-s;\10, \si  )$ do depend on $\si$. 
\item[ii)]  We can view the result of Corollary 8.3 as a new construction of the Kubota-Leopoldt $\ell$-adic $L$-functions.
\end{enumerate}

\medskip

\section{$\ell$-adic  functions associated to measure $K_1(-1)$}

\smallskip 

In this section we identify $\ell$-adic functions
\[
 \Zc _0^\beta (1-s;-1,\si )
\]
constructed with an aid of the measure $K_1(-1)$. 
Let $\varphi$ be a path on $\Pbb _{\bar \Qbb }^1\setminus \{ 0,1,\infty \}$ from $\01$ to $-1$ as on the picture.

\[
 \;
\]

\[
 \;
\]

\[
 \;
\]

$$ {\rm Picture \;6}$$

Let us set
\[
 \delta :=\varphi \cdot x^{\frac{1}{2}}\; .
\]

\medskip

\noindent{\bf Proposition 9.1.} We have
\[
 l(-1)_\delta =0,\;\;li_1(-1)=l_1(-1)_\delta =\kappa (2),
\]
where $\kappa (2)$ is a Kummer character associated to $2$,
\begin{equation}\label{eq:350}
 li_k(-1 )_\delta=l_k(-1)_\delta ={\frac{1-2^{k-1}}{2^{k -1}}}l_k(\10 )_p
\end{equation}
for $k>1$ ($p$ is the standard path from $\01$ to $\10$).

\medskip

\noindent
{\bf Proof.}  The path $\delta$ is chosen so that $l(-1)_\delta =l(\10 )_p=0$. The formula \eqref{eq:350} then follows from the distribution relation
\[
 2^{k-1}\big( l_k(\10 )+l_k(-1)_\delta \big)=l_k(\10 )_p\; ,
\]
whose detailed proof can be found in \cite{NW3}. \hpb

\medskip

From now on we assume that $\ell$ is an odd prime.
Let $\varphi ^{(n)}$ be the path $\varphi _0:=\varphi$ considered on $V_n=\Pbb ^1_{\bar \Qbb }\setminus (\{0,\infty \}\cup \mu _{\ell ^n})$. Let us set
\[
 \delta _n:=\varphi ^{(n)}\cdot x_n^{1/2}
\]
for $n\in \Nbb$ (the loop $x_n$ around $0$ is as in section 2). 
Observe that the constant family $((-1))_{n\in \Nbb }$ is a compatible family of $\ell ^n$-th roots of $-1$.

\medskip

\noindent
{\bf Lemma 9.2.}  We have 
\[
 (\delta _n)_{n\in \Nbb }\in {\varprojlim} _n \pi (V_n;-1,\01 )\, .
\]

\medskip

\noindent
{\bf Proof.}  Let $f:\Cbb ^\times \to \Cbb ^\times $ be given by $f(z)=z^{\ell}$. Then we have
\[
 f(\de )=f(\varphi \cdot x^{1/2})=f(\varphi )\cdot f(x^{1/2})=\varphi \cdot x^{-\frac{\ell-1}{2}}\cdot x^{\frac{\ell}{2}}=\varphi \cdot x^{1/2}=\de \,.
\]
We can assume that all happens in a small neighbourhood of $0$, 
as the image of the interval $[-1,-\varepsilon ]$ ($\varepsilon >0$ and small) is the interval $[-1,-\varepsilon ^{\ell }]$. \hpb

\medskip

It follows from Proposition 2.2 that for $r>0$ we get measures
\[
 K_r(-1)\;.
\]
Hence it follows from Theorem 2.5 (the polylogarithmic case was already proved in \cite{NW}) that
\begin{equation}\label{eq:360}
 l_k(-1)_\delta =li_k(-1 )_\delta= {\frac{1 }{(k-1)!}}\int _{\Zbb _\ell}x^{k-1}dK_1(-1)\,.
\end{equation}
Finally it follows from Proposition 5.3, point i) or the careful examination of the formula v) of Proposition 5.1 that
\begin{equation}\label{eq:370}
 \int _{\Zbb _\ell}x^{k-1}dK_1(-1)={\frac{1}{1-\ell ^{k-1}}}\int _{\Zbb _\ell^\times }x^{k-1}dK_1(-1)\;. 
\end{equation}

\medskip

\noindent
{\bf Definition 9.3.} Let $0\leq \be <\ell -1$. For $\si \in G_\Qbb$ such that $\chi (\si )^{\ell -1}\neq 1$ we define 
\[
 L^\be (1-s;-1,\si ):=
{\frac{2}{\omega (\chi (\si ))^{\be }[\chi (\si )]^s -1}}\int _{\Zbb _\ell ^\times}[x]^s  x^{-1}  \omega (x)^{\beta} dK_1(-1)_\delta (\si)\;.
\]

\medskip

\noindent{\bf Theorem 9.4.} Let $\si \in G_{\Qbb }$ be such that $\chi (\si )^{\ell -1}\neq 1$. 
\begin{enumerate}
 \item[i)] Let $k\equiv \be $ modulo $\ell -1$. Then we have.
\[
  L^\be (1-k;-1,\si )=
{\frac{2  (1-\ell ^{k-1})\cdot (k-1)!}{\chi (\si )^k-1}}l_k(-1)_\delta =
\]
\[
{\frac{2  (1-\ell ^{k-1})\cdot (k-1)!}{\chi (\si )^k-1}}\cdot {\frac{1-2^{k-1}}{2^{k-1}}}\cdot     l_k(\10)_p\;.
\]
\item[ii)] Let $k$ and $\be$ be even and let $k\equiv \be $ modulo $\ell -1$. Then we have
\[
  L^\be (1-k;-1,\si )=(1-\ell^{k-1})\cdot {\frac{1-2^{k-1}}{2^{k-1}}}\cdot {\frac {-B_k}{k}}={\frac{(1-\ell ^{ k-1})  (1-2^{k-1})}{2 ^{k-1}}}  \zeta (1-k)\; .
\]
\end{enumerate}

\medskip

\noindent
{\bf Proof.}  The point i) follows from the formulas \eqref{eq:370}, \eqref{eq:360} and \eqref{eq:350}. The point ii) follows from the point i), the formula  
\eqref{eq:93} and the equality $\zeta (1-k)={\frac{-B_k}{k}}$.   \hpb

\medskip

\noindent{\bf Corollary 9.5.}  Let $\be$ be even and $0\leq \be \leq \ell -3$. Let $\si \in G_\Qbb $ be such that $\chi (\si )^{\ell-1}\neq 1$. 
The function $L^\beta (1-s;-1,\si )$ does not depend on $\si$ and we have
\begin{equation}\label{eq:380}
L^\be (1-s;-1,\si)={\frac{1-2^{-1}  \omega (2)^\beta   [2]^s}{2^{-1}  \omega (2)^\beta   [2]^s}}  L_\ell (1-s,\omega ^\be)\,.
\end{equation}

\medskip

\noindent
{\bf Proof.} Let $\si _1$ and $\si _2$ belonging to $G_\Qbb $ be such that $\chi (\si _1)^{\ell-1}\neq 1\neq \chi (\si _2)^{\ell-1}$. 
Then it follows from Theorem 9.4, ii) that the functions $L^\beta (1-s;-1,\si _1)$ and $L^\beta (1-s;-1,\si _2)$ coinside on the dense subset
\[
 \{k\in \Nbb \mid k\equiv \be \;\;{\rm mod}\;\; \ell-1\}
\]
of $\Zbb _\ell $. Therefore
\[
  L^\beta (1-s;-1,\si _1)=L^\beta (1-s;-1,\si _2)
\]
for any $s\in \Zbb _\ell$.  For $k\in \Nbb$ and $k\equiv \be$ modulo $\ell -1$ it follows from \eqref{lang2} that
\[
 {\frac{1-2^{-1}  \omega (2)^\beta   [2]^k}{2^{-1} \omega (2)^\beta   [2]^k}}  L_\ell (1-k,\omega ^\be)=
{\frac{1-2^{k-1}}{2^{k-1}}}(1-\ell ^{k-1})\zeta (1-k)\,.
\]
Hence   the formula \eqref{eq:380} of the corollary follows from Theorem 2.4, point ii), because the both functions 
coinside on the dense subset $\{k\in \Nbb \mid  k\equiv \be \;{\rm mod}\;\ell-1\}$ of $\Zbb _\ell$. \hpb


\medskip
\section{Hurwitz zeta functions and Dirichlet L-series}

\smallskip 

Let $m$ be a positive integer not divisible by $\ell$. In this section we identify functions
corresponding to measures $K_1(\xi ^i _m)(\si )\mp K_1(\xi ^{m-i} _m)(\si )$.

Let $0\leq \be <\ell -1$ and let $\varepsilon  \in \{1,-1\}$. 
Let us set
\[
  \Zc _0^\beta (1-s;(\xi _m^{-i})+\varepsilon (\xi _m^{i}),\si ):= 
\int _{\Zbb _\ell ^\times}[x]^s   x ^{-1}  \omega (x)^\be d\big( K_1(\xi _m^{-i})(\si )+\varepsilon K_1(\xi _m^{i})(\si )\big)\, .
\]

First we fix paths $\al _i$   from $\01$ to $\xi _m^i$   for $0<i<m$ (see Picture 7).

\[
 \;
\]

\[
 \;
\]

\[
 \;
\]

$${\rm Picture \;7}$$

Let us set 
\[
 \be _i:=\al _i \cdot x^{-\frac{i}{m}}
\]
for $0<i<m$. Observe that then $l(\xi _m^i)_{\be _i}=0$. Hence we have 
\begin{equation}\label{eq:Lambda 444}
 \Lambda _{\be _i}(X,Y )\equiv \sum _{k=1}^\infty l_k(\xi _m^i)_{\be _i} YX^{k-1}\;\;{\rm mod}\;\;\Ic^\prime _2(X,Y)\;. 
\end{equation}

Let $h:\Pbb ^1\setminus \{0,1,\infty \}\to \Pbb ^1\setminus \{0,1,\infty \}$ be given by
\[
 h(\zfk)=1/\zfk \, .
\]

Let us define 
\[
 z:=\Ga ^ {-1}\cdot h  (x)\cdot \Ga\; ,
\]
where $\Ga =\Ga _0$ (see Picture 4).
Then $x\cdot y\cdot z=1$ in $\pi _1(V_0,\01)$.
\medskip

\noindent{\bf Lemma 10.1.} Let $0<i<{\frac{m}{2}}$. Then
\[
 \be _{m-i}=h(\be _i)\cdot \Ga \cdot z^{\frac{i}{m}}\cdot x^{\frac{i}{m}}\,.
\]

\medskip

\noindent
{\bf Proof.} We have
\[
 \be _{m-i}=\al _{m-i}\cdot x^{-{\frac{m-i}{m}}}=\al _ {m-i}\cdot x^{-1}\cdot x^{\frac{i}{m}}=h(\al _i)\cdot \Ga \cdot x^{\frac{i}{m}}=
\]

\[
 h(\al _i\cdot x^{-{\frac{i}{m}}})\cdot h(x^{{\frac{ i}{m}}})\cdot \Ga \cdot x^{\frac{i}{m}}=h(\be _i)\cdot \Ga \cdot z^{\frac{i}{m}}\cdot x^{\frac{i}{m}}\,. 
\]

\hpb


We shall prove the following result.

\medskip

\noindent{\bf Theorem 10.2.} Let $m$ be a positive integer not divisible by $\ell$. We have
\[
 l_k(\xi _m^{-i})_{\be _{m-i}}+(-1)^k l_k(\xi _m^i)_{\be _{ i}}={\frac{1}{k!}}B_k({\frac{i}{m}})\cdot (1-\chi ^k)\,.
\]

\bigskip

To prove Theorem 10.2 we shall need several lemmas.
It follows from Lemma 10.1,   \cite[Lemma 1.0.6]{W1}  and the commuting of $h$ with the action of $G_{\Qbb}$ (see also \cite[formula 10.0.1]{W2}) that
\[
 \ffk _{\be _{m-i}}=\ffk _{h(\be _i) \cdot \Ga \cdot z^{\frac{i}{m}}\cdot x^{\frac{i}{m}}}=
\]
 
\[
  x^{-\frac{i}{m}}\cdot \big( z^{-\frac{i}{m}} \cdot \big( \Ga ^{-1}\cdot h_* (\ffk _{\be _i})\cdot \Ga \cdot \ffk _\Ga \big)\cdot 
  z^{ \frac{i}{m}}\cdot \ffk _{z^{ \frac{i}{m}}}\big)\cdot
  x^{\frac{i}{m}}\cdot \ffk _{x^{\frac{i}{m}}}\, .
\]

We recall that $Z=-\log (\exp X \cdot \exp Y).$ Therefore we get the equality of formal power series
\begin{equation}\label{eq:Lambdaseries}
 \Lambda _{\be _{m-i}}(X,Y)= 
\end{equation}
\[
 e^{{-\frac{i}{m}}X}\cdot \big( e^{{-\frac{i}{m}}Z}\cdot \big(\Lambda _{\be _i}(Z,Y)\cdot \Lambda _\Ga (X,Y)\big) \cdot e^{{\frac{i}{m}}Z}
 \cdot \Lambda _{z^{ \frac{i}{m}}}(X,Y)\big)\cdot e^{{ \frac{i}{m}}X}\cdot e^{{ \frac{i}{m}}(\chi -1)X}\,.
\]
Taking logarithm of both sides of the equality \eqref{eq:Lambdaseries} we get 

\begin{equation}\label{logLambdaseries}
 \log \Lambda _{\be _{m-i}}(X,Y)= \big[ e^{{-\frac{i}{m}}X}\cdot
\end{equation}
\[
  \big( \big(e^{{-\frac{i}{m}}Z}\cdot  [\log \Lambda _{\be _i}(Z,Y)\bigcirc \log \Lambda _\Ga (X,Y) ] \cdot e^{{\frac{i}{m}}Z}\big)\bigcirc
 \log \Lambda _{z^{ \frac{i}{m}}}(X,Y)\big)\cdot e^{{ \frac{i}{m}}X}\big] \bigcirc  {{ \frac{i}{m}}(\chi -1)X}\,.
\]

We shall calculate successive terms of the left hand side of the equality \eqref{logLambdaseries} modulo the ideal $\Ic^\prime _2(X,Y)$.

\medskip

\noindent{\bf Lemma 10.3.} We have
\begin{equation}\label{Lambda zi/m}
  \log \Lambda _{z^{ \frac{i}{m}}}(X,Y)\equiv
\end{equation}
\[
 Y\cdot \Big[ \Big( {\frac{ \exp ( {{ \frac{i}{m}}(1-\chi )X})-\exp ({-{ \frac{i}{m}}\chi   X} ) }{\exp X-1}} 
 +
 {\frac{\chi}{\exp(\chi   X)-1}}\cdot   (e^{-{ \frac{i}{m}}\chi   X}-1) \Big)
\]

\[
 \cdot  {\frac{   {\frac{i}{m}}(1-\chi )X        }{\exp ( {\frac{i}{m}}(1-\chi )  X)-1      }}\Big]
+{\frac{i}{m}}(1-\chi )X\;\;{\rm modulo}\;\;\iii \;.
\]

\medskip

\noindent
{\bf Proof.} We have
$$
\ffk _{z^{ \frac{i}{m}}}(\si )=z^{- \frac{i}{m}}\cdot \si (z^{ \frac{i}{m}})=(x\cdot y)^{\frac{i}{m}}\cdot (\si (x)\cdot \si (y))^{-\frac{i}{m}}\equiv 
(x\cdot y)^{\frac{i}{m}}\cdot (x^{\chi (\si )}\cdot y^{\chi (\si )})^{-\frac{i}{m}}
$$
modulo commutators with two or more $y$'s. Hence we get
\[
 \log \Lambda _{z^{ \frac{i}{m}}}(X,Y)\equiv{ \frac{i}{m}}(X\bigcirc Y)\bigcirc ({-\frac{i}{m}}(\chi   X\bigcirc \chi   Y))\equiv
\]

\[
 \big({ \frac{i}{m}} X+Y\cdot { \frac{{ \frac{i}{m}}  X  }{  \exp X-1}    }\big)\bigcirc 
 \big({ -\frac{i}{m}}\chi   X+Y { \frac{{- \frac{i}{m}}\chi   X  }{  \exp(\chi   X)-1}    }\big)\;\;{\rm mod}\;\; \iii \,.
\]
Applying the formula from Lemma 0.2.1 we get the congruence \eqref{Lambda zi/m} of the lemma. \hpb

\medskip

\noindent{\bf Lemma 10.4.} We have
\[
 \Lambda _\Ga (X,Y)-1\equiv Y\big( {\frac{1}{\exp X-1}}-{\frac{\chi}{\exp (\chi X)-1}}\big)\;\;{\rm mod}\;\;\iii\;.
\]

\medskip

\noindent
{\bf Proof.} Observe that $\Ga =h(p)^{-1}\cdot s\cdot p$. Hence we have \[\ffk _\Ga=\Ga ^{-1}\cdot h_*(\ffk _p^{-1})\cdot \Ga 
                                                                         \cdot p^{-1}\cdot \ffk _s\cdot p\cdot \ffk _p\;.\]  
Therefore after the embedding of $\pi _1(\Pbb ^1 _{\bar \Qbb }\setminus \{0,1,\infty \},\01 ) $
into $\Qbb _{\ell}\{\{X,Y\}\}  $ we get
\[
  \Lambda _\Ga (X,Y)=\Lambda _p(Z,Y)^{-1}\cdot e^{{\frac{1}{2}}(\chi -1)Y}\cdot \Lambda _p(X,Y)\,.
\]
Hence it follows from the congruence \eqref{eq:defpoly} that
\[
 \log  \Lambda _\Ga (X,Y)=(-\log \Lambda _p(Z,Y))\bigcirc ({{\frac{1}{2}}(\chi -1)Y})\bigcirc \log  \Lambda _p(X,Y)\equiv
\]

\[
 \big(\sum _{k=2}^\infty(-1)^kl_k(\10 )_pYX^{k-1}\big)\bigcirc \big({{\frac{1}{2}}(\chi -1)Y}\big)\bigcirc \big(\sum _{k=2}^\infty l_k(\10 )_pYX^{k-1}\big)\equiv
\]

\[
 {{\frac{1}{2}}(\chi -1)Y}+\sum _{k=1}^\infty 2 l_{2k}(\10 )_pYX^{2k-1}\;\;{\rm mod}\;\;\iii\,.
\]
In \cite{W7} we have shown that
\[
 l_{2k}(\10 )_p={\frac{B_{2k}}{2\cdot (2k)!}}(1-\chi ^{2k})
\]
(see also \cite[Proposition 5.13]{NW2}). Therefore we get
\[
  \log  \Lambda _\Ga (X,Y)\equiv \sum _{k=1}^\infty {\frac{B_k}{k!}}(1-\chi ^k)YX^{k-1}\;\;{\rm mod}\;\;\iii \,.
\]
It follows from the definition of the Bernoulli numbers that the right hand side of the last congruence is equal
\[
 Y\big( {\frac{1}{\exp X-1}} -{\frac{1}{X}}\big) -Y\big(   {\frac{\chi}{\exp (\chi X)-1}}-{\frac{1}{X}}\big) = 
 Y\big( {\frac{1}{\exp X-1}} -  {\frac{\chi}{\exp (\chi X)-1}}\big) \,.
\]
It is clear that $\Lambda _\Ga (X,Y)-1\equiv \log \Lambda _\Ga (X,Y)$ modulo $\iii$. Hence the lemme follows. \hpb

\medskip

\noindent
{\bf Proof of Theorem 10.2.}  Let us set
\[
 A_i(X):=\sum _{k=1}^\infty l_k(\xi _m^i)_{\be _i}X^{k-1}\;.
 \]
Observe that
\[
 \log \Lambda _{\be _i}(Z,Y)\bigcirc \log \Lambda _\Ga (X,Y)\equiv Y\Big( A_i(-X)+{\frac{1}{e^ X-1}}-{\frac{\chi}{ e^{\chi  X}-1}}\Big) \;{\rm mod} \;\iii
\]
and
\begin{equation}\label{zBi}
 e^{-{\frac{i}{m}}Z}\big( \log \Lambda _{\be _i}(Z,Y)\bigcirc \log \Lambda _\Ga (X,Y)\big)  e^{{\frac{i}{m}}Z}\equiv
\end{equation}
\[
 Y\Big( A_i(-X)+{\frac{1}{\exp X-1}}-{\frac{\chi}{\exp(\chi  X)-1}}\Big)    e^{-{\frac{i}{m}}X}\;\;{\mod}\;\;\iii\;.
\]
Let us denote by \[ S(X) \]
the formal power series in the square bracket of the congruence \eqref{Lambda zi/m} of Lemma 10.3,
i.e. we have
\[
  \log \Lambda _{z^{ \frac{i}{m}}}(X,Y)\equiv Y  S(X)+{\frac{i}{m}}(1-\chi )X\;\;{\rm mod}\;\;\iii\;.
\]
It follows from the congruences \eqref{zBi} and \eqref{Lambda zi/m} and Lemma 0.2.1 that
\begin{equation}\label{eq:*}
 e^{-{\frac{i}{m}}X}\cdot \Big( ( e^{-{\frac{i}{m}}Z}\cdot (\log  \Lambda _{\be _i}(Z,Y)\bigcirc \log \Lambda _\Ga (X,Y) )\cdot e^{{\frac{i}{m}}Z})\bigcirc
 \log \Lambda _{z^{ \frac{i}{m}}}(X,Y)\Big) \cdot e^{{\frac{i}{m}}X}\equiv
\end{equation}
$$Y\cdot \Big( (A_i(-X)+{\frac{1}{\exp X-1}}-{\frac{\chi}{\exp (\chi X)-1}})\cdot  e^{-{\frac{i}{m}}X}\cdot$$
\[ 
{\frac{\exp ( {\frac{i}{m}} (1-\chi )X)\cdot  {\frac{i}{m}} (1-\chi )X}{\exp ( {\frac{i}{m}} (1-\chi )X)-1}}+S(X)\Big)\cdot e^{{\frac{i}{m}}X}+{\frac{i}{m}}(1-\chi )X\;\;{\rm mod}\;\;\iii\,.
\]
Following the equality \eqref{logLambdaseries} it rests to calculate the $\bigcirc$-product of the right hand side of \eqref{eq:*} with 
${\frac{i}{m}}(1-\chi )X$. Using once more Lemma 0.2.1 we get
\begin{equation}\label{eq:final}
 \log \Lambda _{\be_{m-i}}(X,Y)\equiv Y  \Big(A_i(-X)+{\frac{\exp ( {\frac{i}{m}} X) }{\exp X-1}}-{\frac{\chi \exp ( {\frac{i}{m}} \chi X) }{\exp (\chi X)-1}}\Big)\;\;{\rm mod}\;\;\iii\;. 
\end{equation}
We recall that the Bernoulli polynomials $B_k(t)$ are defined by the generating function
\[
 {\frac{X  \exp(tX)}{\exp X-1}}=\sum _{k=0}^\infty{\frac{B_k(t)}{k!}}  X^k\;.
\]
Therefore finally we get the following congruence
\begin{equation}
 Y  \big( \sum _{k=1}^\infty l_k(\xi _m^{m-i})_{\be _{m-i}}X^{k-1}\big)\equiv
\end{equation}
\[
 Y  \Big( \sum _{k=1}^\infty (-1)^{k-1}l_k(\xi _m^{i})_{\be _{i}}X^{k-1}+\sum _{k=1}^\infty {\frac{B_k({\frac{i}{m}    }  )}{k!}}\cdot (1-\chi ^k)X^{k-1}\Big)\;.
\]
Comparing the coefficients we get
\[
 l_k(\xi _m^{-i})_{\be _{m-i}}+(-1)^kl_k(\xi _m^i)_{\be _i}={\frac{B_k({\frac{i}{m}    }  )}{k!}}  (1-\chi ^k)\;.
\]
\hpb

\medskip

\noindent
{\bf Proposition 10.5.} Let $m$ be a positive integer not divisible by $\ell$. We have
\begin{equation}\label{eq:prop10.6}
 {\frac{1}{1-\chi ^k}}\int _{\Zbb _\ell}x^{k-1}d(K_1(\xi _m^{-i})+(-1)^{k}K_1(\xi _m^i))={\frac{B_k({\frac{i}{m}    }  )}{k }}
\end{equation}
for $0<i<m$ and $k\geq 1$.

\medskip

\noindent
{\bf Proof.} For $0<i<{\frac{m}{2}}$ the proposition follows immediately from Theorem 2.5 (see also \cite[Proposition 3]{NW}). 
If ${\frac{m}{2}}<i<m$ then we use the equality $B_k(1-X)=(-1)^kB_k(X)$ (see \cite[page 41]{H}).  \hpb

\medskip

We recall here the definition of Hurwitz zeta functions. Let $0<x\leq 1$. Then one defines
\[
 \zeta (s,x):=\sum _{n=0}^\infty (n+x)^{-s} 
\]
(see \cite[page 41]{H}). The function $\zeta (s,x)$ can be continued beyond the region $\Re (s)>1$. One shows that
\[
 \zeta (1-n,x)=-{\frac{B_n(x)}{n}}
\]
for all $n>0$ (see \cite[Section 2.3, Theorem 1]{H}). We shall construct $\ell$-adic non-Archimedean analogues of the Hurwitz zeta functions using measures 
$K_1(\xi _m^{-\al })\pm  K_1(\xi _m^\al )$.

\medskip
Let $\al ={\frac{a}{b}}$ be a rational number and let $a$ and $b$ be integers. We assume that $b$ and $m$ are relatively prime. 
Then we define the integer $\langle \al \rangle$ by the conditions $0\leq \langle \al \rangle <m$ and 
$\langle \al \rangle \equiv \al $ modulo  $m$.

\medskip

\noindent
{\bf Proposition 10.6.} Let $m$ be a positive integer not divisible by $\ell$. Let $a$ be the order of $\ell$ in $(\Zbb /m\Zbb )^\times$. Let $0<\al <m$ be such that $(\al ,m)=1$. Then we have
\begin{equation}\label{eq:Prop10.7}
{\frac{1}{1-\chi ^k}} \int _{\Zbb _\ell ^\times}x^{k-1}d\Big( K_1(\xi _m^{-\al \ell ^{-p}})+(-1)^k K_1(\xi _m^{\al \ell ^{-p}})\Big)= 
\end{equation}
\[
 {\frac{1}{k}}\Big( B_k( {\frac{\langle \al \ell ^{-p}\rangle }{m}})-\ell ^{k-1}  B_k(  {\frac{\langle \al \ell ^{-p-1}\rangle }{m}})\Big)
\]

for $p=0,1,\ldots a-1$.

\medskip

\noindent
{\bf Proof.} Observe that
\[
 \int _{\Zbb _\ell }x^{k-1}dK_1(\xi _m^\al )=\sum _{i=0}^{a-1}{\frac{\ell^{(k-1)i}}{1-\ell ^{(k-1)a}}}\int _{\Zbb _\ell ^\times }x^{k-1}K_1(\xi _m^{\al {\ell ^{-i}}})
\]
by Proposition 5.3. Hence it follows from Proposition 10.5 that
\[
 (1-\ell ^{(k-1)a})   {\frac{1}{k}} B_k( {\frac{\langle \al \ell ^{-j}\rangle }{m}})=
\] 
 
\[
 {\frac{1}{1-\chi ^k}}\sum _{i=0}^{a-1}\ell ^{(k-1)i}
 \int _{\Zbb _\ell ^\times }
 x^{k-1}d\Big( K_1(\xi _m^{-\al \ell ^{-j-i}})+(-1)^k K_1(\xi _m^{\al \ell ^{ {-j-i}}})\Big) 
\]
for $j=0,1,\ldots ,a-1$. Multiplying the $(p+1)$th equation by $\ell ^{k-1}$ and next subtracting from the $p$th equation and dividing by $(1-\ell ^{(k-1)a})$ we get
the equalities \eqref{eq:Prop10.7} of the proposition. \hpb

\medskip

\noindent
{\bf Remark 10.6.1} A similar formula as the right hand side of equalities  \eqref{eq:Prop10.7} appears in \cite[Theorem 1]{Sh}.
\medskip

Let $\varepsilon \in \{1,-1\}$. We define
\[
 L^ \be (1-s;(\xi _m^{- i})+\varepsilon (\xi _m^ {i}),\si ):= {\frac{1}{\omega (\chi (\si ))^ \be [\chi (\si )]^ s-1}}\cdot 
 \Zc ^\be _0(1-s;(\xi _m^{-i})+\varepsilon (\xi _m^i),\si )=
\]

\[
 {\frac{1}{\omega (\chi (\si ))^ \be [\chi (\si )]^ s-1}}\int _{\Zbb _{\ell}^ \times }[x]^ s   x^ {-1}  \omega (x)^ \be 
 d\big(K_1(\xi _m^{-i })(\si )+\varepsilon  K_1(\xi _m^i )(\si )\big)\,.
\]

\medskip

\noindent
{\bf Proposition 10.7.} Let $0\leq \be <\ell -1$ and let $\si \in G_{\Qbb}$ be such that $\chi (\si )^{\ell-1}\neq 1$. Then for $k\equiv \be$ modulo $\ell -1$ we have
\[
  L^ \be (1-k;(\xi _m^{- i})+(-1)^\be(\xi _m^ { i}),\si )=
  {\frac{1}{k}}\Big( B_k( {\frac{\langle i\rangle }{m}})-\ell ^{k-1}  B_k(  {\frac{\langle i \ell ^{-1}\rangle }{m}})\Big)\,.
\]

\medskip

\noindent
{\bf Proof.} The proposition follows immediately from Proposition 10.6. \hpb

\medskip

\noindent
{\bf Corollary 10.8.} Let $\si $ and $\si _1$ be such that $\chi (\si )^{\ell -1}\neq 1$ and $\chi (\si _1)^{\ell -1}\neq 1$. Then we have
\[
 L^ \be (1-s;(\xi _m^{- i})+(-1)^\be(\xi _m^ { i}),\si )\,=\,L^ \be (1-s;(\xi _m^{- i})+(-1)^\be(\xi _m^ { i}),\si _1)\,.
\]

\medskip

\noindent
{\bf Proof.} Both functions take the same values at the dense subset of $\Zbb _\ell$, hence they are equal. \hpb

\medskip

\noindent
{\bf Remark 10.8.1.}  Notice that  the function 
$L^ \be (1-s;(\xi _m^{- i})+(-1)^\be(\xi _m^ { i}),\si )$ does not depend on the choice of $\si$ such that $\chi (\si )^{\ell -1}\neq 1$. 
This function is then an  $\ell$-adic 
non-Archimedean analogues of the Hurwitz zeta function $\zeta (s,{\frac{i}{m}})$.

\medskip

Let $\psi :(\Zbb /q\Zbb )^\times \to \bar \Qbb ^\times$ be a primitive Dirichlet character. The L-series attached to $\psi$ is defined by
\[
 L(s,\psi )=\sum _{n=1}^\infty {\frac{\psi (n)}{n ^s}} 
\]
for $\Re (s)>1$. Then one shows that 
\[
 L(s, \psi )=\sum _{a=1}^q \psi (a)  q^{-s}  \zeta (s,{\frac{a}{q}})
\]
and for $n>1$ one has 
\[
 L(1-n,\psi )=-{\frac{1}{n}}q^{n-1}\sum _{a=1}^q \psi (a)\cdot B_n({\frac{a}{q}})
\]
(see \cite[Chapter 4, page 31 and Theorem 4.2.]{Wa}).

\medskip

Having $\ell$-adic non-Archimedean Hurwitz zeta functions we shall define $\ell$-adic Dirichlet L-series.
Let $m$ be a positiveinteger not divisible by $\ell$. Let
\[
 \psi :(\Zbb/m\Zbb)^\times \to \bar \Qbb _\ell ^\times
\]
be a primitive Dirichlet character. Let $0\leq \beta <\ell -1$ and let $\varepsilon \in \{1,-1\}$. 
Let $\si \in G_{\bar \Qbb }$ be such that $\chi (\si )^{\ell -1}\neq 1$. We define
\[
 L_\ell ^\beta (1-s;\psi ,\varepsilon,\si ):=-\omega (m)^\beta [m]^sm^{-1}\sum _{a=1}^m\psi (a)   L^\beta (1-s;(\xi _m^{-a})+\varepsilon (\xi ^a _m),\si )\,.
\]
\medskip
\medskip

\noindent
{\bf Proposition 10.9.}
\begin{enumerate}
 \item[i)] The function $ L_\ell ^\beta (1-s;\psi ,(-1)^ \beta ,\si ) $ does not depend on a choice of $\si \in G_\Qbb $.
 \item[ii)] For $k\equiv \beta $ modulo $\ell -1$ we have 
 \[
   L_\ell ^\beta (1-s;\psi ,(-1)^ \beta ,\si )=(1-\psi (\ell)  \ell ^ {k-1})  L(1-k,\psi )\,.
 \]
\end{enumerate}
\medskip

\noindent
{\bf Proof.} We calculate
\[
 L_\ell ^\beta (1-s;\psi ,(-1)^ \beta ,\si )=-\omega (m)^ \beta [m]^ km^ {-1}\sum _{a=1}^m\psi (a)   L^\beta (1-k;(\xi _m^{-a})+(-1)^ \beta (\xi ^a _m),\si)= 
\]
\[
 -m^ {k-1}\sum _{a=1}^ m\psi (a){\frac{1}{k}}\big(B_k({\frac{a}{m}})-\ell ^ {k-1}B_k({\frac{\langle a\ell ^ {-1}\rangle }{m}})\big)=
\]
\[
 -m ^ {k-1}{\frac{1}{k}}\big(\sum _{a=1}^ m\psi (a)B_k({\frac{a}{m}} )-\ell ^ {k-1}\sum _{a=1}^ m\psi (\ell )\psi (a\ell ^ {-1})B_k({\frac{\langle a\ell ^ {-1}\rangle }{m}})\big)=
\]
\[
 L(1-k,\psi )-\ell ^ {k-1}\psi (\ell)L(1-k,\psi )=(1-\psi (\ell)\ell ^ {k-1})  L(1-k,\psi ).
\]
Hence we have proved the point ii). The first statement is now clear.  \hpb

\medskip
\noindent
{\bf Remark 10.10.} If $\varepsilon \neq (-1)^\beta$ then the functions $ L_\ell ^\beta (1-s;\psi ,\varepsilon,\si )$ do depend on $\si \in G_{\Qbb }$. 
We think that the measure 
\[
 \sum _{a=1}^m\psi (a)\big( K_1(\xi _m^{-a})(\si )+\varepsilon K_1(\xi _m^a)(\si )\big)
\]
can be called $\ell$-adic Dirichlet L-series of the character $\psi$. The measure $K_1(\10 )$ is then $\ell$-adic zeta function.
In fact these measures can be considered as measures on $\hat \Zbb $ not only on $\Zbb _\ell$ (see \cite{NW} and also \cite{W8}).


\medskip
\section{$\ell$-adic L-functions of $\Zb [1/m]$}

\smallskip 

The functions $L_\ell (1-s;-1,\si )$ considered in Section 9 can be view as the $\ell$-adic L-function of $\Zb [1/2]$. 
Let $p_1,p_2,\ldots ,p_r$ be different prime numbers. Below we propose to define an $\ell$-adic L-functions of  $\Zb [1/m]$.

\medskip

\noindent
{\bf Lemma 11.1.} Let $p_1,p_2,\ldots ,p_r$ be different prime numbers. Let $m= p_1  p_2\ldots p_r$.  Then we have
\[
 \sum _{i=1, (i,m)=1}^{m-1}B_k\big(\langle {\frac{i}{m}}\rangle \big)=\Big( \prod _{j=1}^r{\frac {(1-p_j^{k-1})}{p_j^{k-1}}}\Big)  B_k\,.
\]

\medskip

\noindent
{\bf Proof.} The distribution formula for Bernoulli polynomials implies the equality
\begin{equation}\label{eq:BER1}
m^{k-1}\big( \sum _{i=0}^{m -1} B_k({\frac {i}{m}})\big)=B_k\,.
\end{equation}
Let $P:=\{p_1,p_2,\ldots ,p_r\}$. If $A=\{p_{a_1},\ldots ,p_{a_s}\}$ is a subset of $P$ we set 
\[
N_A:={\frac{p_1  p_2\ldots p_r}{p_{a_1}  p_{a_2}\ldots  p_{a_s}}}\, . 
\]
Then we can write the equality \eqref{eq:BER1} in the form
\begin{equation}\label{eq:BER2}
m^{k-1}\big( \sum _{i=0,\, (i,m)=1}^{m -1} B_k({\frac {i}{m}})
+\sum _{\emptyset \neq A\subset P}\sum _{i=0,\,(i,N_A)=1}^{N_A-1}B_k({\frac{i}{N_A}})\big)=B_k\,.
\end{equation}
The equality \eqref{eq:BER1} implies immediately the formula of the lemma for $r=1$. Let us suppose that the formula of the lemma is true for all $q<r$. 
Then we get from the equality \eqref{eq:BER2} the following equality
\begin{equation}\label{eq:BER3}
m^{k-1}  \sum _{i=0,\, (i,m)=1}^{m -1} B_k({\frac {i}{m}})
+\sum _{\emptyset \neq A\subset P}\big(\prod _{p\in A}p^{k-1}\big)    \big( \prod _{p\in P\setminus A}(1-p^{k-1})\big) B_k=B_k\,.
\end{equation} 
Let {\bf r} be the set $\{1,2,\ldots ,r\}$. Let us write the Taylor formula for the polynomial $X_1X_2\ldots   X_r$ at the point $(p_1^{k-1},\ldots ,p_r^{k-1})$.
We get 
\[
 X_1X_2 \ldots   X_r =p_1^{k-1} p_2^{k-1} \ldots   p_r^{k-1}+ \sum _{\emptyset \neq B\subsetneqq {\bf r}}
 \big(\prod _{j\in {\bf r}\setminus B}p_j^{k-1}\big)   \big( \prod _{i\in B}(X_i-p_i^{k-1})\big)+
 \prod _{i=1}^r(X_i-p_i^{k-1}).
\]
Setting $(X_1,\ldots ,X_r)=(1,\ldots ,1)$ we get
\begin{equation}\label{eq:T1}
 \prod _{i=1}^r(1-p_i^{k-1})+\sum _{B\subsetneqq {\bf r}}\big(\prod _{j\in {\bf r}\setminus B}p_j^{k-1}\big) \cdot \big( \prod _{i\in B}(1-p_i^{k-1})\big)=1\,.
\end{equation}
Comparing the equalities \eqref{eq:BER3} and \eqref{eq:T1} we get the equality of the lemma. \hpb

\medskip

For $0\leq \be <\ell -1$ and $\si \in G_\Qbb$ such that $\chi (\si )^{\ell -1}\neq 1$ we define 
\[
 L^\be (1-s,\Zbb [{\frac{1}{m}}],\si ):=
\]
\[
 {\frac{2}{\omega  (\chi (\si ))^\be   [\chi (\si )]^\be -1}}\cdot 
 \int _{\Zbb _\ell ^\times }[x]^s   x^{-1}  \omega (x)^\be d\big(\sum _{i=1,\, (i,m)=1}^{m-1}K_1(\xi _m^{-i})(\si )\big)\,.
\]
Let us assume that $\beta$ is even and $k\equiv \beta$ modulo $\ell -1$. From the very definition of the function 
$L^\be (1-s,\Zbb [{\frac{1}{m}}],\si)$ we have
\[
 L^\be (1-k,\Zbb [{\frac{1}{m}}],\si )= \sum _{i=1,\, (i,m)=1}^ mL^ \beta (1-k;(\xi _m^ {-i})+(-1)^ \beta (\xi _m^ i),\si )\,.
\]
Hence it follows from Proposition 10.7 and Lemma 11.1 that
\[
 L^\be (1-k,\Zbb [{\frac{1}{m}}],\si )=(-1)^ r {\frac{1}{k}}  (1-\ell ^ {k-1})B_k  \prod _{j=1}^ r (p_j   p_j^ {-k}-1)\,.
\]
Hence it follows from the equality \eqref{lang2} that

\begin{equation}\label{last}
 L^\be (1-k,\Zbb [{\frac{1}{m}}],\si )=L_\ell (1-k,\omega ^ \beta)  \prod _{j=1}^ r(p_j [p_j]^ {-k}  \omega (p_j)^ {-\beta }-1)\,.
\end{equation}

\medskip

\noindent
{\bf Proposition 11.2.} Let $p_1,p_2,\ldots ,p_r$ be different prime numbers and let $m=p_1\cdot p_2\ldots p_r$. Let $\beta $ be even and let $0\leq \beta <\ell -1$. 
Let $\si \in G_\Qbb $ be such that $\chi (\si )^ {\ell -1}\neq 1$. Then we have 
\[
 L^\be (1-s,\Zbb [{\frac{1}{m}}],\si)=\prod _{j=1}^ r\big( p_j  [p_j]^ {-s}\cdot \omega (p_j)^ {-\beta }-1\big)  L_\ell (1-s,\omega ^ \beta )\,.
\]
\medskip

\noindent
{\bf Proof.} The proposition follows immediately from the equality \eqref{last}. \hpb

\medskip

\noindent
{\bf Proposition 11.3.} Let $p$ be a prime number. Then we have
\[
 \int _{\Zbb _\ell}[x]  x^{-1}  \omega (x)  d\big( \sum _{i=1}^{p-1}K_1(\xi _p^{-1})(\si )\big)=l(p)(\si )\,.
\]
\medskip

\noindent
{\bf Proof.}  The integral is equal $\sum _{i=1}^{p-1}l_1(\xi _p^{-i})(\si )=l(p)(\si )$. \hpb

\medskip

Notice that 
$$\int _{\Zbb _\ell ^\times }d K_1(\xi _p^j)(\si )= l(1-\xi _p^j)(\si ) -l(1-\xi _p^{j\ell ^{-1}})(\si )\, ,$$ hence 
$\int _{\Zbb _\ell ^\times }d(\sum _{i=1}^{p-1} K_1(\xi _p^{-i})(\si ))=0$.

In view of Proposition 11.2 and 11.3 we can consider the measure $\sum _{i=1}^{p-1} K_1(\xi _p^{-i})$ as an $\ell$-adic zeta function of the ring 
$\Zbb [{\frac{1}{p}}]$. However if $m$ is a product of $r$
different prime numbers with $r>1$ then the integral $\int _{\Zbb _\ell}d\big( \sum _{i=1,\, (i,m)=1}^{p-1} K_1(\xi _m^{-i})(\si )\big)=0$, 
but dim$_{\Qbb _\ell}H^1(\Zbb [{\frac{1}{m}}];\Qbb _\ell (1))=r$. We can replace the measure $\sum _{i=1,\, (i,m)=1}^{p-1} K_1(\xi _m^{-i}) $ by the measure 
$\sum _{i=1 }^{p-1} K_1(\xi _m^{-i}) $. Then $\int _{\Zbb _\ell }d\big( \sum _{i=1 }^{p-1} K_1(\xi _m^{-i})(\si )\big) =l(m)(\si )$ and
$$
{\frac{2}{\omega (\chi (\si ))^\beta   [\chi (\si )]^\beta -1}}  \int _{\Zbb _\ell ^\times }[x]^s 
  x^{-1}  \omega (x)^\beta d\big( \sum _{i=1 }^{p-1} K_1(\xi _m^{-i})(\si )\big)
=
$$
\[
\big( m  [m]^{-s}  \omega (m)^{-\beta }-1\big)  L_\ell (1-s,\omega ^\beta )
\]
if $\beta$ is even and $\chi(\si )^{\ell-1}\neq 1$. 
We do not know which choice is better if any.

\bigskip

\noindent
The results of this paper were presented in the international meeting on polylogarithms in June 2012 in Nice 
and in the poster session of the Iwasawa 2012 conference in Heidelberg.

\bigskip

\noindent
{\bf Acknowledgment} These research were started in January 2011 during our visit in Max-Planck-Institut f\"ur Mathematik in Bonn. 
We would like to thank very much MPI for support.

\vglue 2cm

\vglue 1cm

\noindent Universit\'e de Nice-Sophia Antipolis

\noindent D\'epartement de Math\'ematiques

\noindent Laboratoire Jean Alexandre Dieudonn\'e

\noindent U.R.A. au C.N.R.S., N$^{\rm  o}$ 168

\noindent Parc Valrose -- B.P. N$^{\rm  o}$ 71

\noindent 06108 Nice Cedex 2, France

\smallskip

\noindent {\it E-mail address} wojtkow@math.unice.fr

\noindent {\it Fax number} 04 93 51 79 74

\medskip

\end{document}